\documentclass[iicol]{snjnl}
\usepackage{natbib}
\setcitestyle{authoryear,round}
\usepackage{graphicx}%
\usepackage{multirow}%
\usepackage{amsmath,amssymb,amsfonts}%
\usepackage{amsthm}%
\usepackage{mathrsfs}%
\usepackage[title]{appendix}%
\usepackage{xcolor}%
\usepackage{textcomp}%
\usepackage{manyfoot}%
\usepackage{booktabs}%
\usepackage{algorithm}%
\usepackage{algorithmicx}%
\usepackage{algpseudocode}%
\usepackage{listings}%
\usepackage{subcaption}
\usepackage{caption}

\theoremstyle{thmstyleone}%
\theoremstyle{thmstyletwo}%

\theoremstyle{thmstylethree}%

\begin{document}


\title[Article Title]{Stress and frequency optimization of prismatic sandwich beams with joints: Performance improvements through topology optimization}


\author[1]{\fnm{Shengyu} \sur{Yan}}\email{jasonsy.yan@ubc.ca}

\author*[2]{\fnm{Jasmin} \sur{Jelovica}}\email{jasmin.jelovica@ubc.ca}

\affil[1]{\orgdiv{Department of Mechanical Engineering}, \orgname{The University of British Columbia}, \orgaddress{\street{6250 Applied Science Ln}, \city{Vancouver}, \postcode{V6T 1Z4}, \state{BC}, \country{Canada}}. ORCID: 0000-0002-9077-9649}

\affil[2]{\orgdiv{Departments of Mechanical and Civil Engineering}, \orgname{The University of British Columbia}, \orgaddress{\street{6250 Applied Science Ln}, \city{Vancouver}, \postcode{V6T 1Z4}, \state{BC}, \country{Canada}}. ORCID: 0000-0002-8396-941X}


\abstract{Prismatic sandwich panels fabricated from metals offer a compelling alternative to more traditional panels across diverse industries, primarily due to their superior strength-to-weight ratio. Although several core types were proposed in the past, further improvements in performance could be achieved by devising the topology of the core through a topology optimization framework, which is explored in this article for the first time. Another novelty is the inclusion of joints between the sandwich beams and its surroundings in the analysis and optimization. Stress is minimized under uniform pressure loading on the beams and natural frequency maximized using the Method of Moving Asymptotes. The results are compared with X-core, Y-core, corrugated-core and web-core sandwich beams, a few conventional prismatic sandwich types, which are optimized using a prominent global evolutionary algorithm. Manufacturing requirements are considered through practical limitations on the design variables. It is shown that structures produced by topology optimization outperform the conventional sandwich beams by up to 44\% at intermediate to high mass levels, where volume fraction is between 0.2 and 0.4, but often through increased topological complexity. The new core topologies bear a certain resemblance with the conventional core types, underscoring engineering ingenuity that went into their development over the years. The topology of the optimized joints differs from the conventional joint. The results also show some limitations of the topology optimization framework, for example that it does not offer better performing beams for volume fractions below 0.2. }




\maketitle

\section{Introduction}\label{sec1}

Advanced structural solutions promote performance and sustainability to a higher level in many industries, i.e. automotive, construction, marine, aeronautical and aerospace. Sandwich structure is one type of advanced structure since it provides enhanced energy absorption properties, high stiffness-to-weight ratio, excellent ballistic resistance performance, and good thermal and acoustic isolation properties \citep{birman2018review}. Over the past several decades, sandwich panels have attracted substantial research interest. Sandwich panels consist of two continuous faces (top face and bottom face) and a core. The core can be fabricated from a range of materials and adopt diverse forms, e.g. honeycomb, unidirectional profile, or even foam. Noor et al. offer a comprehensive review of honeycomb and foam-filled sandwich panels, and overview of computational techniques for their structural analysis \citep{noor1996computational}. In the last two decades, research has intensified around sandwich panels with prismatic cores, which are relatively simple to manufacture from metals and suitable for large-scale applications \citep{roland2004advanced},\citep{bright2007new}. Several types of cores have been proposed, e.g., web-core, corrugated-core, X-core, Y-core, etc. Earlier studies have focused on the mechanical response of idealized panels with various cores types, e.g., \citep{chong1982free}, \citep{sakata1982natural}, \citep{sakata1988vibrations}, \citep{li1992forced}, \citep{dalaei1996natural}, \citep{lok2001free}, \citep{lok2000free},\citep{lok2001bending}, whereas more recent investigations considered influence of production on response and strength, see e.g. \citep{romanoff2007laser}, \citep{jelovica2012influence}, \citep{yan2021buckling}. Zhang and Khalifa highlight the high stiffness-to-weight ratio and high impact strength of web-core sandwich panels \citep{zhang2017numerical}, \citep{khalifa2018performance}. Yan and Jelovica studied the buckling and free vibration of web-core sandwich panels \citep{yan2021buckling}. Bending stress analysis of corrugated-core, honeycomb-core and X-core sandwich panels have been studied by He et al. \citeyearpar{he2012precise}. Kiliccaslan et al. revealed the behavior of corrugated core sandwiches under quasi-static and dynamic loadings \citep{kiliccaslan2016single}.

The optimization of metallic sandwich structures continues to be a significant focus in research. Rathbun et al. conducted a thorough strength optimization with bending loads, revealing minor variations in the weights of optimized panels \citep{rathbun2005strength}. Yuan et al. develop failure maps and applied sequential linear programming (SLP) to optimize truss-core panels subjected to thermal loading \citep{yuan2016failure}. Similarly, Sadollah et al. employed the water cycle algorithm (WCA) for sizing optimization of prismatic sandwich panels \citep{sadollah2013sizing}. Further advancements in optimization methods were presented by Cai et al., who used the non-dominated sorting genetic algorithm (NSGA-II) to identify the optimal design for corrugated core sandwich panels under air blast loading conditions \citep{cai2021multi}. Commercial finite element software is sometimes used in such research. As an example, Tauhiduzzaman utilized the ABAQUS shape optimization module to optimize four different types of sandwich cores \citep{tauhiduzzaman2016design}. Sizing optimization has emerged as a robust tool in the quest for optimal design of prismatic sandwich structures, particularly under diverse boundary conditions and load cases. This approach has been significantly enhanced by the use of genetic algorithms. This paper also applies genetic algorithm for sizing optimization, which is then compared with topology optimization. Although certain amount of literature exists on the performance and characteristics of prismatic sandwich structures, information regarding the joining of such structures to the surrounding structures is extremely limited. Joints must be meticulously designed to ensure adequate fatigue life and load-carrying capacity. Several joints for prismatic sandwich panels have been analyzed in \citep{niklas2008search}, with stress concentration factors determined using commercial FE code and various joints proposed based on engineering ingenuity.The study was part of the EU-funded SANDWICH project, where a more comprehensive report was made on the joints and their stress concentration factors by Kujala and Ehlers \citep{SANDWICH}. Yan and Jelovica presented a study on the joints connecting sandwich panels with surrounding girders using topology optimization \citep{yan2022topology}, but the joints were considered in isolation. It is worth noting that numerous studies have been published on composite adhesive joints between steel and composite plates, with examples found in \citep{cen2017mechanical}. However, metallic sandwich panels are intended for use in metallic structures, and as such should be joined to metallic girders or bulkheads, necessitating substantially different joints. 

Topology optimization (TO) is a prominent method for structural optimization, which could propose structures with better performance than the traditional ones. In recent years, great efforts have been made in this area and led to a series of approaches, such as the density method \citep{bendsoe1989optimal}, level set method \citep{allaire2002level}, evolutionary approaches \citep{xie1993simple}, topological derivative \citep{sokolowski1999topological}, phase field method \citep{bourdin2003design} and machine learning \citep{lei2019machine}. Many scholars contributed to the optimization of composite sandwich structures \citep{sun2017topological}, lattice core sandwich \citep{xiao2021design} and honeycomb sandwich structures \citep{sun2017topological}. However, a study on topology optimization of prismatic metallic sandwich structures has not been reported. The authors have presented a study on the joint connecting sandwich beams with the surrounding girder structures using topology optimization \citep{yan2022topology}, but the joint was considered in isolation. However, stiffness of the sandwich beam could affect the topology of the joints, and thus they need to be optimized together. 

In this paper, we perform topology optimization of prismatic sandwich beams, with and without joints to the surrounding structure, all made from metal. Stress is minimized \citep{deng2021efficient} and natural frequency is maximized \citep{du2007topological} using the Method of Moving Asymptotes (MMA) \citep{svanberg1987method}. The results are compared to panels with conventional cores: corrugated-, web-, X- and Y-core sandwich panels. Sizing optimization (SO) is performed on the conventional beams, incorporating joint proposed by Kujala and Ehlers \citep{SANDWICH}. Conventional beams are compared amongst themselves and to the TO results. The capabilities of the TO framework to converge to meaningful results are assessed. The findings contribute to the development of sandwich structures and their joints. 

\section{Methodology}\label{sec2}
\subsection{Topology optimization (TO)}\label{subsec1}
The process of topology optimization is commonly paired with the finite element method (FEM) to discretize the design domain \citep{sigmund2013topology}. In this paper, the density approach is utilized \citep{bendsoe1989optimal}, defining design variables as the relative density value of each element; here, a value of 0 denotes void material and 1 indicates solid material. The values that lie between 0 and 1 are penalized using the amended Solid Isotropic Material with Penalization (SIMP) method \citep{bendsoe2003topology}, effectively penalizing the elastic modulus as
\begin{equation}
    E_i=E_i(x_i)=E_{min}+x_i^{pl}(E_0-E_{min})
\end{equation}
where \(E_i\) is the elastic modulus of \textit{i}th element, \(E_{min}\) is the minimum elastic modulus value (\(10^{-9}\) in this paper), \(E_0\) is the material elastic modulus, \(pl\) is the penalization power $(pl>1)$, and \(x_i\) is the relative density of the \textit{i}th element.

Current topology optimization is stated as follows:
\begin{equation}
\begin{aligned}
    \text{Find}& \quad \mathbf{x} = \begin{bmatrix} x_1,x_2 \dots x_n \end{bmatrix}, \quad x_i \in (0,1]\\
    \text{Minimize}& \quad \sigma_{PN} \quad or \quad \Bar{\lambda}\\
    \text{Subject to}& \quad \mathbf{v(x)}=\mathbf{x}^T\mathbf{v}-\Bar{v}<=0
    \end{aligned}
\end{equation}
where \textbf{x} is the relative density vector, \(\sigma_{PN}\) is the p-norm von Mises stress value, \(\bar{\lambda}\) is the average value of the first \textit{n} natural frequencies, \textbf{v} is the element volume vector and \(\bar{v}\) is the volume fraction (VF).

The robust and computationally efficient Method of Moving Asymptotes (MMA) \citep{svanberg1987method} is adopted as the optimization algorithm due to its capacity to handle diverse objectives and constraints. As a gradient-based optimization algorithm, MMA necessitates the first-order derivatives of all the objective and constraint functions, which are commonly referred to as `sensitivities' within the realm of TO. These sensitivities function as the directional guides in the MMA search process, with the sensitivities from the three most recent iterations factored into each subsequent iteration.

To effectively minimize the peak stress within the structure, an aggregate stress measure becomes imperative. Consequently, the p-norm von Mises stress is employed for this purpose \citep{deng2021efficient}. A rudimentary minimization function proves unsuitable due to its non-differentiable nature; however, the p-norm effectively makes the objective function smoother. The calculation of the p-norm von Mises stress can be achieved as follows:
\begin{equation}
    \sigma_{PN}=(\sum_{i=1}^n = \hat{\sigma}_{vm,i})^{1/p}
\end{equation}
where \(p\) is the norm order and \(\hat{\sigma}_{vm,i}\) is the panelized von Mises stress of the \textit{i}th element, which can be obtained by
\begin{equation}
\begin{split}
    \sigma_{vm,i} = &(\frac{1}{2}[(\sigma_{xx,i} - \sigma_{yy,i})^2 + (\sigma_{yy,i} - \sigma_{zz,i})^2 \\&+ (\sigma_{zz,i} - \sigma_{xx,i})^2 + \\&6(\sigma_{xy,i}^2 + \sigma_{yz,i}^2 + \sigma_{xz,i}^2)])^{1/2}
\end{split}
\end{equation}
and
\begin{equation}
    \hat{\sigma}_{vm,i}=x_i^q\sigma_{vm,i}
\end{equation}
where \(\sigma_{vm,i}\) is the von Mises stress of the \textit{i}th element, and \(q=1.5\) is the stress relaxation factor presented by Holmberg et al.\citep{holmberg2013stress}.

The sensitivity of the p-norm von Mises stress could be obtained by using the chain rule, and the resultant equation is \citep{deng2021efficient}
\begin{equation}
    \frac{\partial \sigma_{PN}}{\partial x_j} = \sum_{i=1}^n \frac{\partial \sigma_{PN}}{\partial \hat{\sigma}_{vm,i}} \left[ \left( \frac{\partial \hat{\sigma}_{vm,i}}{\partial \mathbf{\hat{\sigma_i}}} \right)^T \frac{\partial (x_i^q\mathbf{\sigma}_i)}{\partial x_j} \right]
\end{equation}
where \(\mathbf{\hat{\sigma_i}}\) is the relaxed stress component \(i\).

To optimize the first natural frequency of the structure, a weighted average scheme is implemented to formulate the objective function. This approach considers multiple natural frequencies by using a weight vector in order to relieve the influence of mode shift \citep{du2007topological}.
\begin{equation}
    \Bar{\lambda}=\mathbf{w\lambda}
\end{equation}
where \(\mathbf{w}\) is the row vector of weights and \(\lambda\) is the column vector of eigenfrequencies. The derivative of the \(j\)th eigenfrequency about design variable \(x_i\) can be calculated as
\begin{equation}
    \frac{\partial \lambda_j}{\partial x_i}=\mathbf{\varphi}_j^T(\mathbf{K}_i'-\lambda_j\mathbf{M}_i')\mathbf{\varphi}_j
\end{equation}
where \(\mathbf{\varphi}_j\) is the \(j\)th eigenvector, \(\mathbf{K}_i\) and \(\mathbf{M}_i\) are the penalized stress and mass matrix of the \textit{i}th element, and they can be obtained by
\begin{equation}
    \mathbf{K}_i=x_i^{pl}\mathbf{K}_{i0}
\end{equation}
where \(\mathbf{K}_{i0}\) is the original stiffness matrix of the \textit{i}th elemement.

From a preliminary investigation, it is found that when the mass and stiffness matrices are penalized using the identical factor, \(pl\), localized vibrations tend to occur in regions with lower relative density. This phenomenon, termed `localized modes' by Pedersen, is triggered by a small mass-to-stiffness ratio when the relative density is low \citep{du2007topological}. To provide a more balanced penalization and circumvent the occurrence of spurious eigenmodes, a novel scheme \citep{du2007topological} for penalizing the two matrices has been adopted in the frequency-based TO. This scheme is implemented as follows:
\begin{equation}
    \mathbf{M}_i=
    \begin{cases}
        x_i\mathbf{M}_{i0},& x_i>0.1\\
        x_i^{ql}\mathbf{M}_{i0},\quad &x_i<=0.1
    \end{cases}
\end{equation}
where $ql$ is the penalization power for the mass matrix and $ql>1$.
As a result, their derivatives are
\begin{equation}
\begin{split}
        &\mathbf{K}_i'=pl*x_i^{pl-1}\mathbf{K}_{i0}\\
        &\mathbf{M}_i'=
    \begin{cases}
        \mathbf{M}_{i0}, &x_i>0.1\\
        qlx_i^{ql-1}\mathbf{M}_{i0},\quad &x_i<=0.1
    \end{cases}
\end{split}
\end{equation}

The effectiveness and efficiency of TO can be influenced by a series of control parameters and factors, i.e. element size, filter radius, stiffness penalization power, mass penalization power, etc, which use the values in the research of \cite{du2007topological} and \cite{deng2021efficient}. A sensitivity study on the element size will be discussed in Section \ref{subsec1}.

\subsection{Sizing optimization (SO)}\label{subsec2}
Sizing optimization is executed utilizing the Non-Dominated Sorting Genetic Algorithm II (NSGA-II) \citep{deb2002fast}, which is a benchmark multi-objective optimization algorithm. To avoid any potential bias in the results, each optimization run is initiated from randomly generated designs. For every optimization scenario, two runs are executed. The resultant non-dominated designs from each run are amalgamated and filtered via non-dominated sorting to obtain the streamlined front. As per our preliminary study, there is no requirement for a local search akin to the method used in \citep{yan2021buckling} as the global optimization algorithm demonstrates complete success in identifying all the Pareto optimal designs for our case study.

The variables considered in the optimization process include the thicknesses of the face plate, thickness of the plates in the core, height of the core, and spacing of the plates in the core. The objectives taken into account are maximizing the 1st natural frequency, minimizing the maximum von Mises stress, and minimizing the weight (outside view area density of a sandwich beam). This optimization process does not include explicit constraints. The variables are treated as continuous values, allowing them to vary smoothly and to have more permutations in the given bounds. Typical cross-over and mutation probabilities are used \citep{deb2002fast}. The algorithm is executed for 100 generations, with each generation consisting of a population of 120.

\begin{figure}[h]
    \centering
    \includegraphics[width=.49\textwidth]{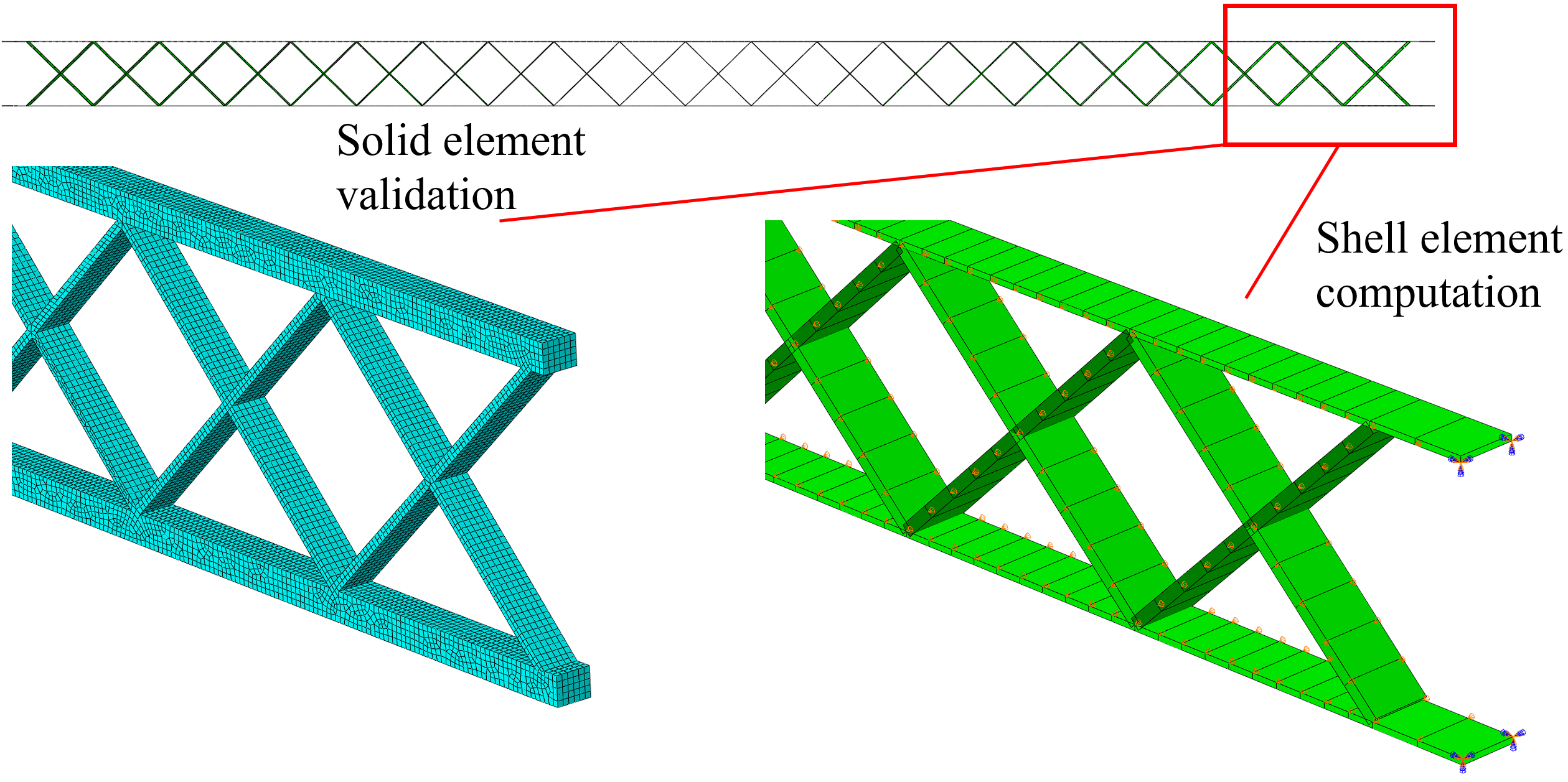}
    \caption{Abaqus shell and solid models used in sizing optimization, shown on X-core sandwich beam as example.}
    \label{XcoreAbaqus}
\end{figure}

The stress and frequency of the prismatic sandwich beams are obtained through FEA models using 4-node shell elements (S4R) with reduced integration for computational efficiency in Abaqus R2018. An example of the shell model and solid model of the X-core beam is shown in Figure \ref{XcoreAbaqus}. Mesh density is uniformly defined within each representative volume element (RVE) of each core type. A range from 8 to 20 elements per unit length has been tested in our preliminary study, and 10 to 12 is found to have the best efficiency and accuracy, thus 10 is selected for the further study. 

\subsection{Validation}\label{subsec3}
The validation is conducted using Abaqus R2018. Two kinds of higher-order models have been built: (1) 3D model using 8-node solid elements, C3D8R, for the SO results and (2) 2D model using 4-node shell elements, S4R, for the TO results. Results from randomly selected cases have been validated, and all the cases show a maximum difference in stress or frequency values within 3\%.

\section{Case studies}\label{sec3}
\subsection{General data}
\begin{figure}[hb]
    \centering
    \includegraphics[width=.49\textwidth]{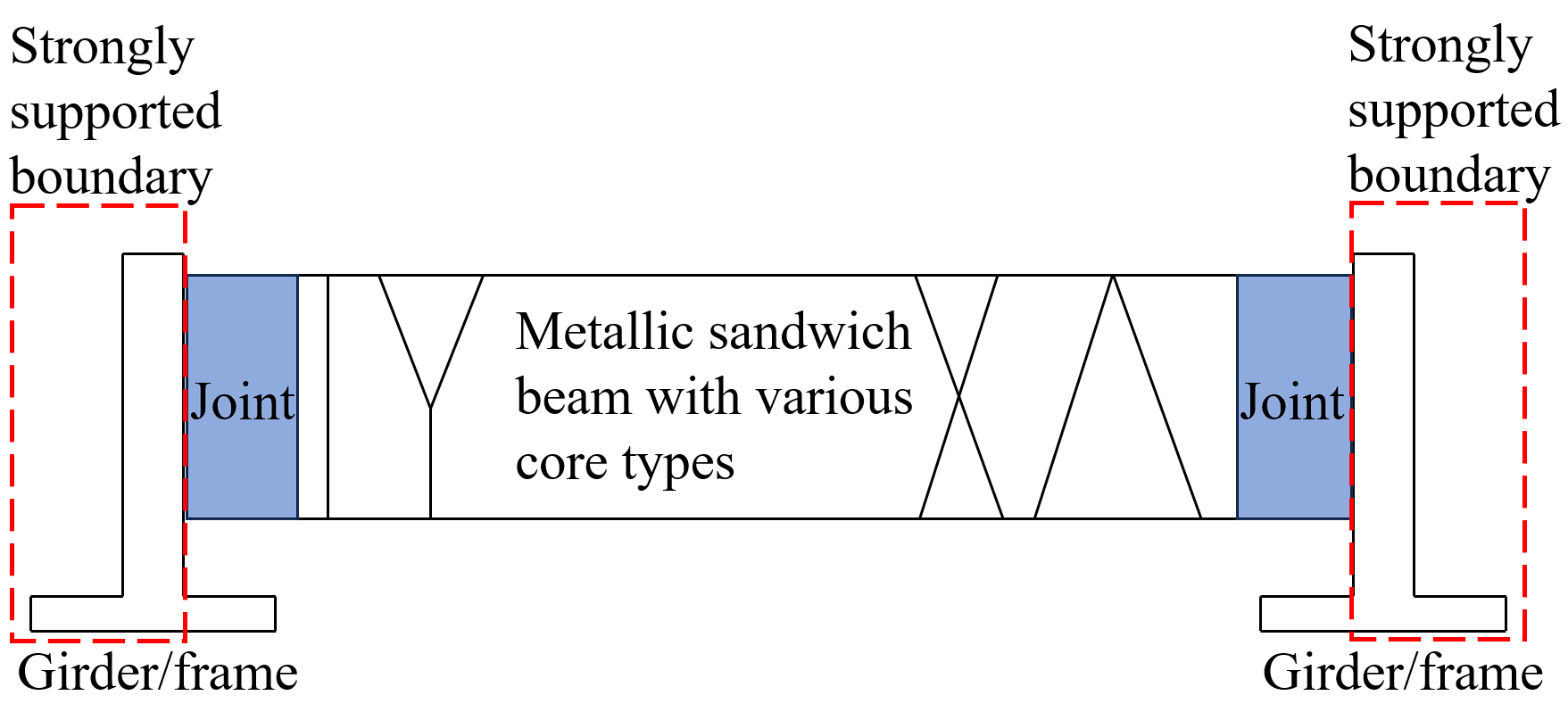}
    \caption{Metallic prismatic sandwich beams connected to the surrounding structures.}
    \label{deck}
\end{figure}

We consider steel beams with material density of 7850 \(kg/m^3\). The beam is used here to represent the cross-section of a panel typically spanning the distance between two nearby girders or frames in ships, bridges, and similar structures, as shown in Figure \ref{deck}. The distance between girders typically starts from 1 m, which is thus the length of the beam used in this study. In order to simplify the complexity of this analysis, a linear elastic material law has been implemented. Two sets of objective functions are evaluated and optimized separately - (1) minimum mass and highest 1st natural frequency, and (2) minimum mass and minimum highest von Mises stress. The assumed boundary conditions for all scenarios are identical: both ends of the beam are fully clamped, and a distributed pressure load of 50 kPa is applied to the upper surface in the stress simulation. The finite element analysis (FEA) in SO is executed via the Abaqus R2018 software, using shell elements (S4R). On the other hand, TO employs MATLAB and utilizes linear quadrilateral elements (S4). The influence of the choice of FEA solver on the results is elaborated in Section \ref{sec4}.

\begin{figure}[h]
    \centering
    \includegraphics[width=.45\textwidth]{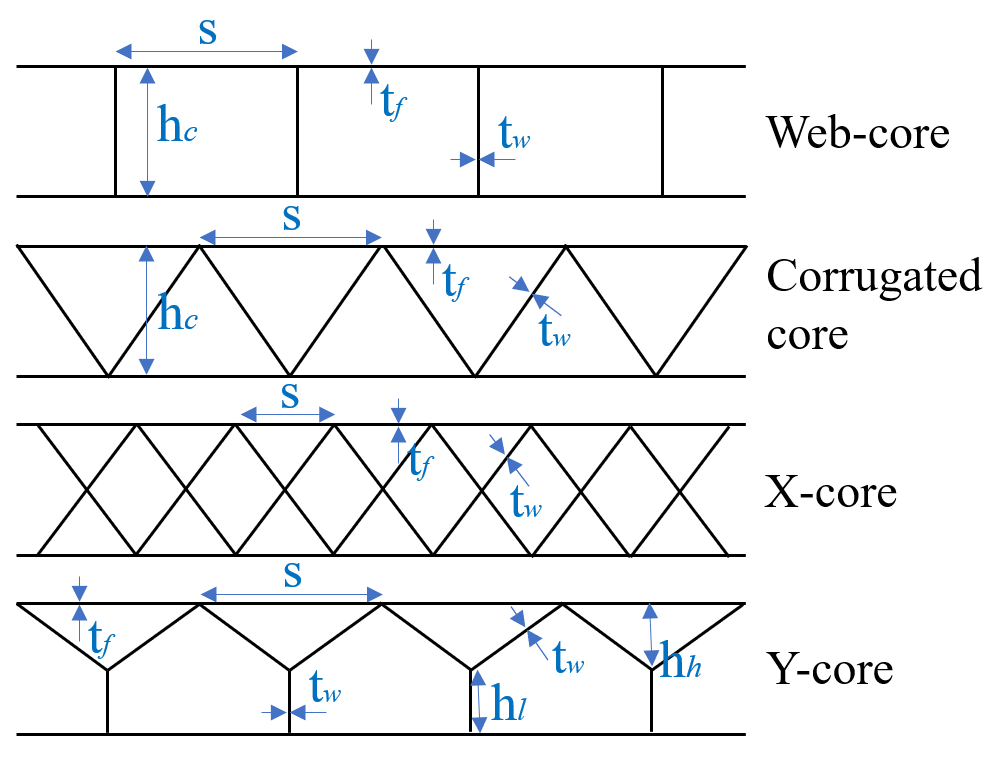}
    \caption{The four kinds of metallic prismatic sandwich beam core types used in SO with variable dimensions.}
    \label{coretypes}
\end{figure}

Four types of cores are studied in the SO: corrugated core, web-core, X-core, and Y-core as shown in Figure \ref{coretypes}. The design variables are annotated on the figure as well. Variable bounds are given in Table \ref{tab1}, following range of values used in \cite{jelovica2012influence}. Lower limits for the thickness of plates is set to the smallest values that could be used in case of e.g. laser welding. 

\begin{table}[h]
\captionsetup{justification=raggedright,singlelinecheck=false, width=\linewidth}
\caption{Variable bounds used in SO.}\label{tab1}
\centering
\begin{tabular}{lcc}
\toprule
Variables & Min [mm]  & Max [mm] \\
\midrule
$t_f$ & 1.0 & 4.0 \\
$t_w$ & 2.5 & 5.0 \\
$h_c$ & 10 & 37 \\
$h_h^*$ & 1 & 36 \\
$h_l^*$ & 1 & 36 \\
s & 30 & 150 \\
\bottomrule
\end{tabular}
\raggedright\small *only used in the Y-core type, where \(h_h+h_l \leq 37\)
\end{table}

\subsection{Sandwich beams without joints}
\begin{figure}[h]
    \centering
    \includegraphics[width=0.315\textwidth]{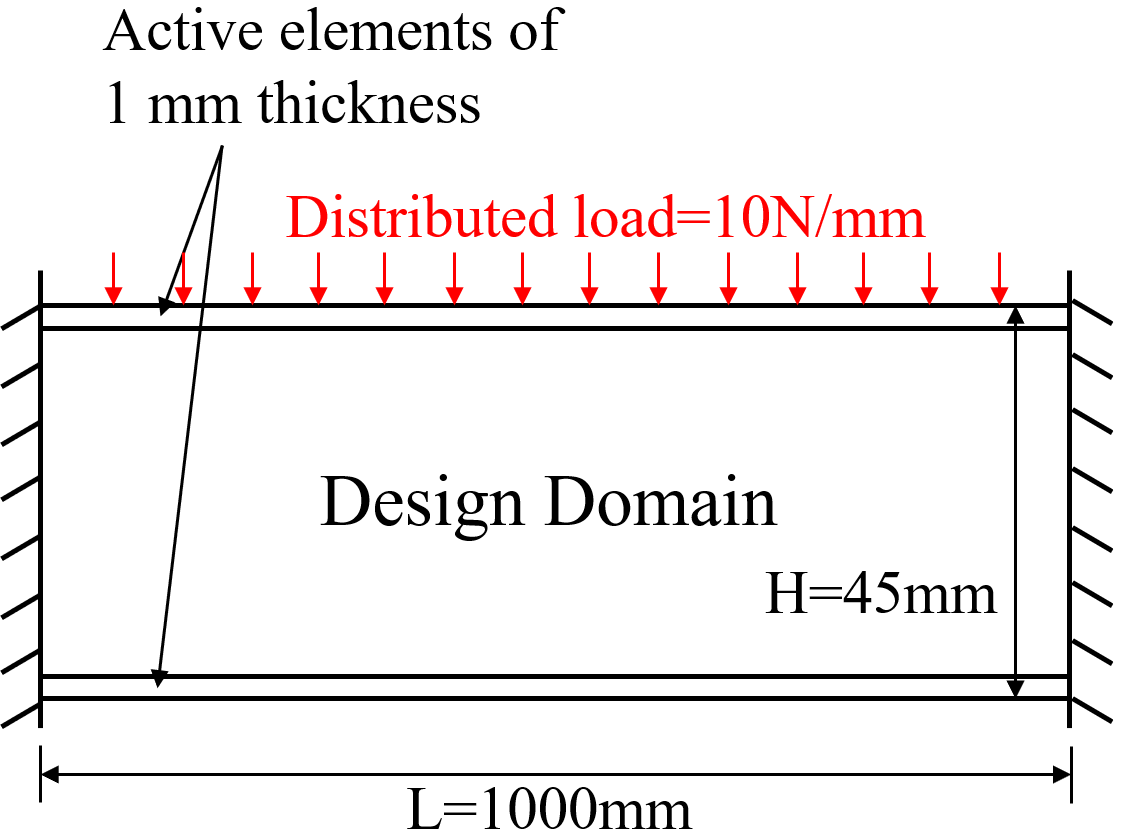}
    \caption{Boundary conditions of the beam without joints in TO.}
    \label{NojointTOBC}
\end{figure}

The scope of our investigation encompasses two main case studies. The first addresses sandwich beams without joints, seeking to explore the beam topologies that offer the best performance, independent of fixture considerations. The design domain in this case is illustrated in Figure \ref{NojointTOBC}.

\subsection{Sandwich beams with joints}
The second case study considers sandwich beams with joints. Topology optimization was performed on sandwich beam joints for a few loading cases \citep{yan2022topology}. However, the joint was optimized independently, and no beam or panel structure is considered adjacent to it. The question is whether the optimal configuration or topology will alter when the design domain expands to include the surrounds. Consequently, we have undertaken stress minimization optimization using both topology optimization and sizing optimization in this research, in an endeavor to gain deeper insights into this problem. Frequency optimization is excluded in this case, as our preliminary study suggest that the joint does not considerably influence the global eigenpairs, due to the negligible alterations in the relationship between mass and stiffness.

\begin{figure}[h]
    \centering
    \includegraphics[width=0.45\textwidth]{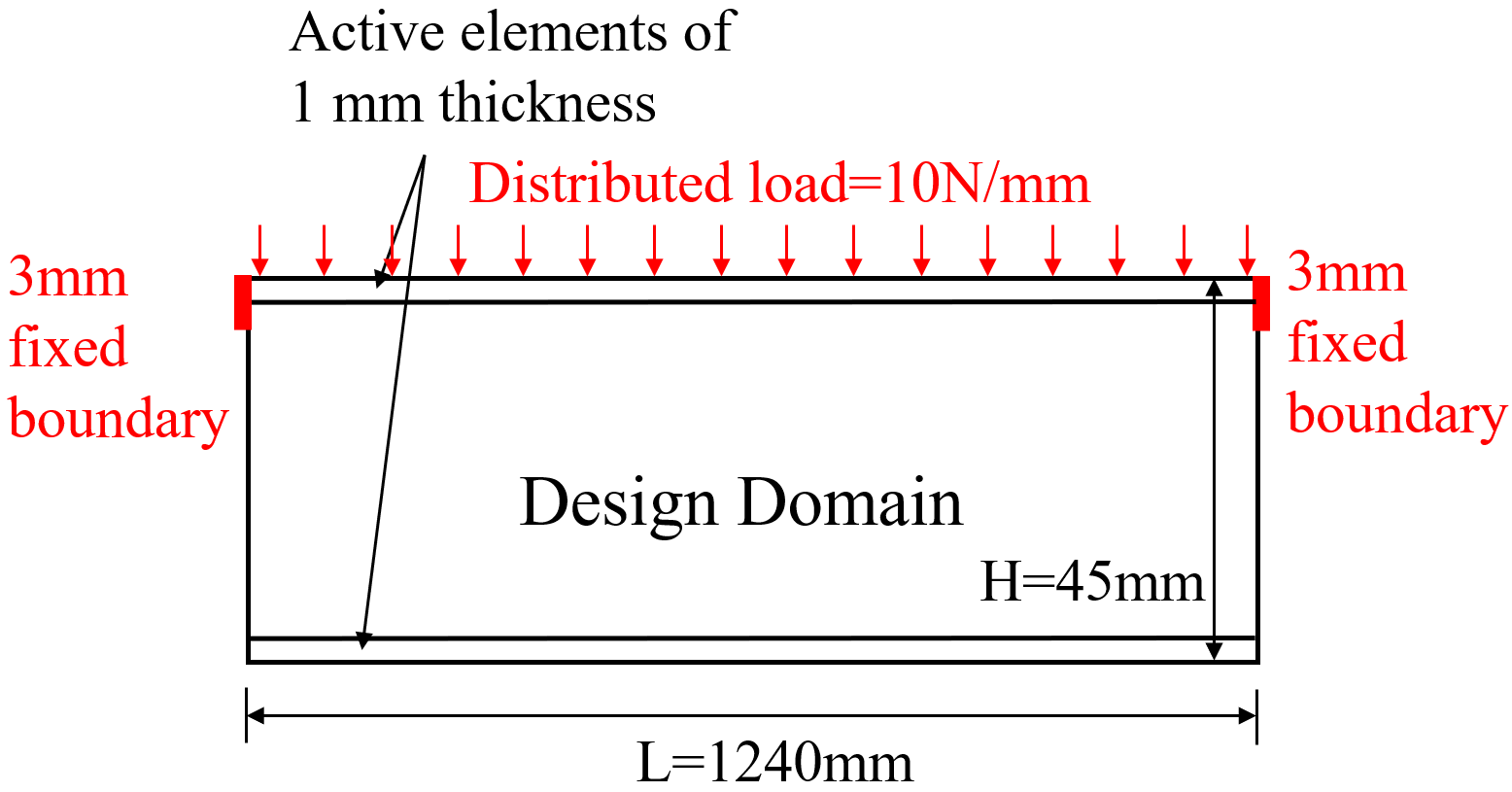}
    \caption{Boundary conditions of the beam with joints in TO.}
    \label{JointTOBC}
\end{figure}

\begin{figure}[h]
    \centering
    \includegraphics[width=0.45\textwidth]{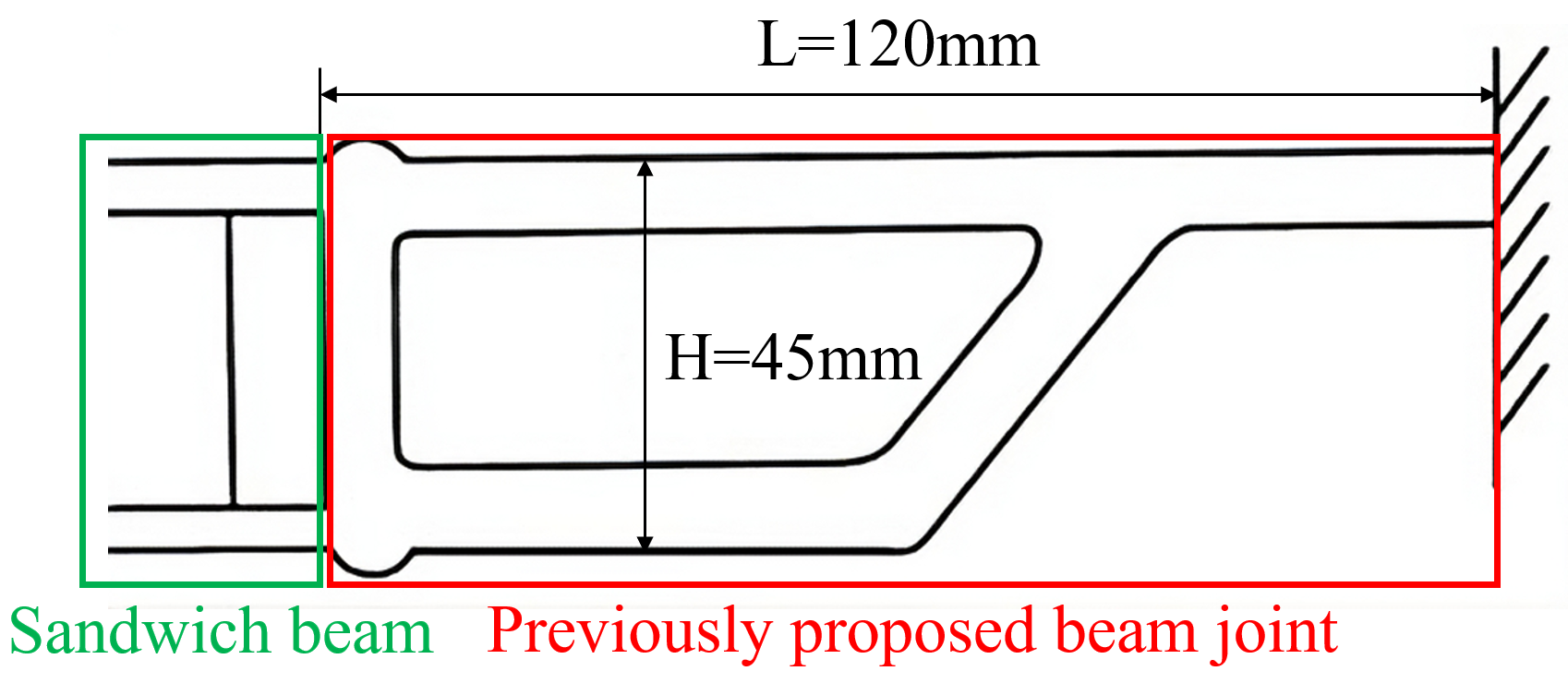}
    \caption{Previously proposed joint connecting the metallic sandwich beam to girder/frame.}
    \label{convetionaljoint}
\end{figure}

\begin{figure}[h]
    \centering
    \includegraphics[width=0.49\textwidth]{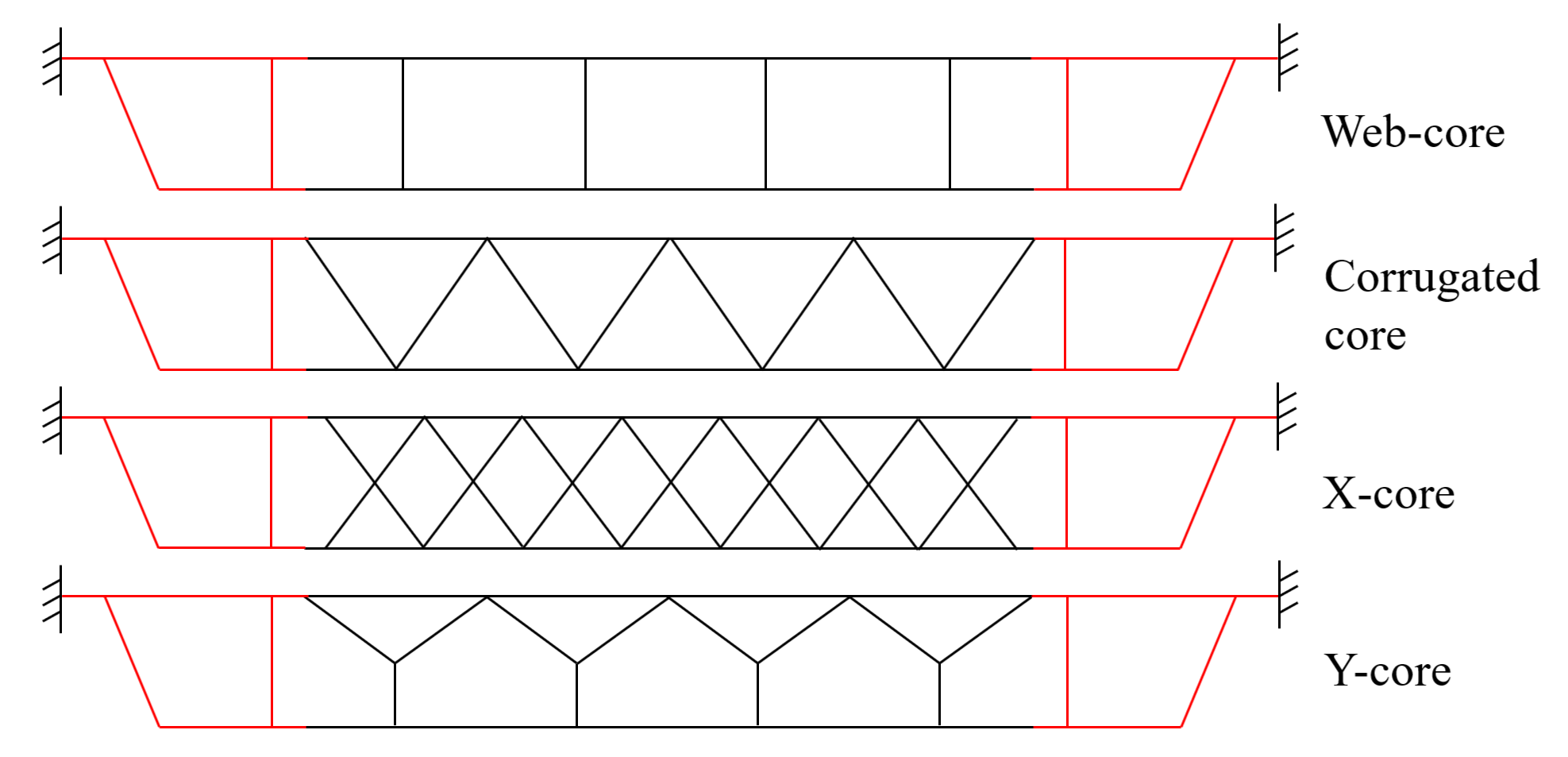}
    \caption{Prismatic sandwich beams integrated with previously proposed joints in SO.}
    \label{core type with joints}
\end{figure}

The design domain for topology optimization in this new context is depicted in Figure \ref{JointTOBC}. Simultaneously, sizing optimization uses previously proposed joint for the sandwich beams, as shown in Figure \ref{convetionaljoint}. The beams with the four different core types are integrated with the joints on both ends as shown in Figure \ref{core type with joints}. In SO, the thickness of the joint is also introduced as a new design variable, $t_j$, and $1mm \leq t_j \leq 3mm$. The new design domain is larger for both TO and SO in this case, as two joint sections have been added to the one-meter-long beam.

\section{Results and discussion}\label{sec4}
\subsection{Sandwich beams without joints}\label{sec:nojoint}
\subsubsection{Stress optimization}\label{sec:stressnojoint}
\begin{figure}[h]
    \centering
    \includegraphics[width=0.49\textwidth]{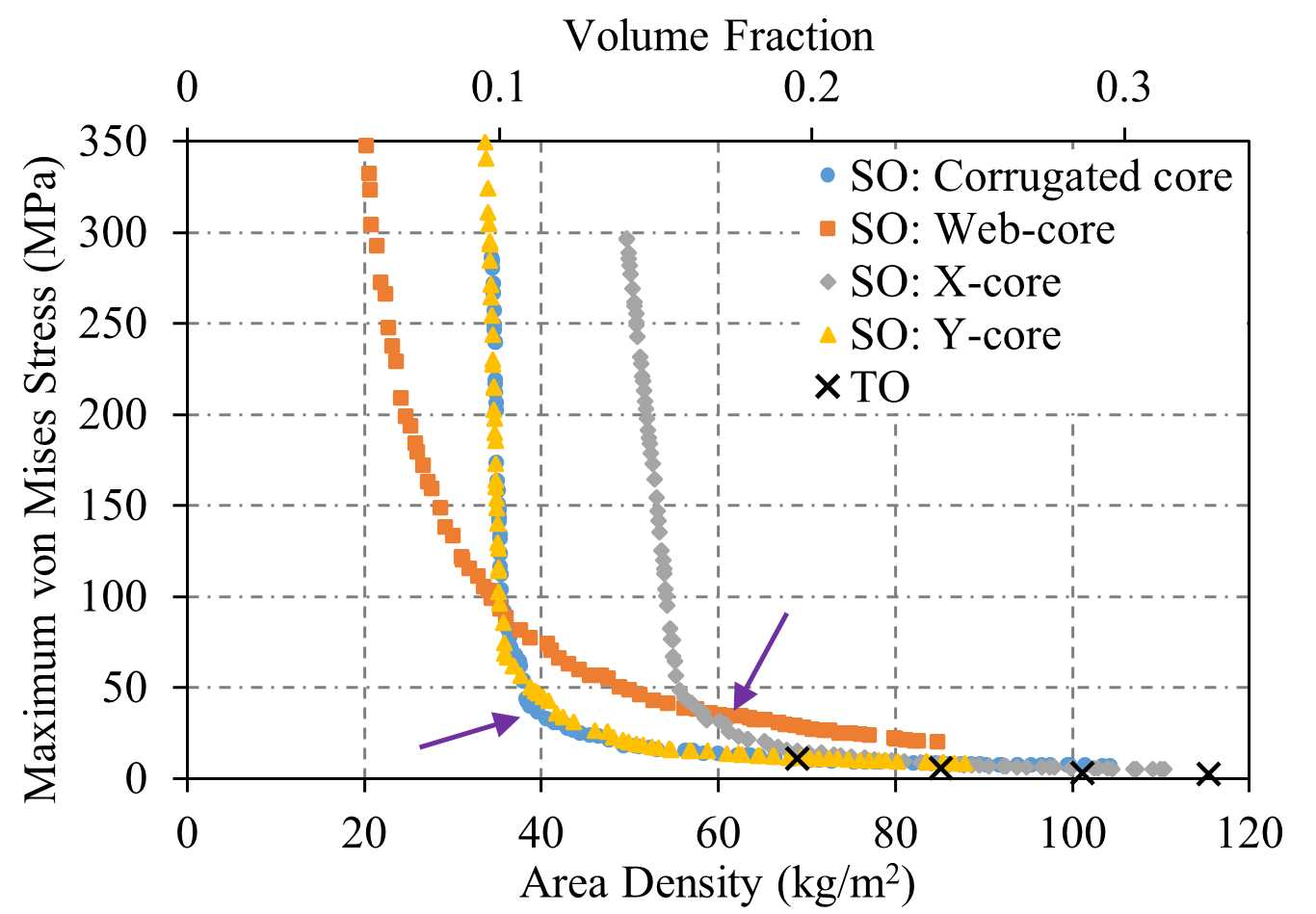}
    \caption{Stress-based optimization results for beams without joints.}
    \label{noJstress}
\end{figure}

Pareto fronts of the optimized sandwich beams are illustrated in Figure \ref{noJstress}. It can be seen that NSGA-II successfully discerns broad and smooth range of solutions for the traditional prismatic beams. TO produces very competitive results. Among SO solutions, web-core sandwich beams reach the lowest area densities, but with relatively high stress. At higher area densities, above 35 kg/m2, the other three sandwich types are better at minimizing the stress. The performance of the corrugated-core and Y-core beams is relatively similar, since the Y-core beams in the Pareto front mostly possess minimum value of the lower core height (thh), resulting in a core shape that resembles corrugated core. X-core sandwich beams are generally the heaviest, and are inferior to corrugated-core beams in terms of the maximum stress, except above 85 kg/m2 (for VF above 0.25). 

\begin{figure}[h]
    \centering
    \includegraphics[width=0.49\textwidth]{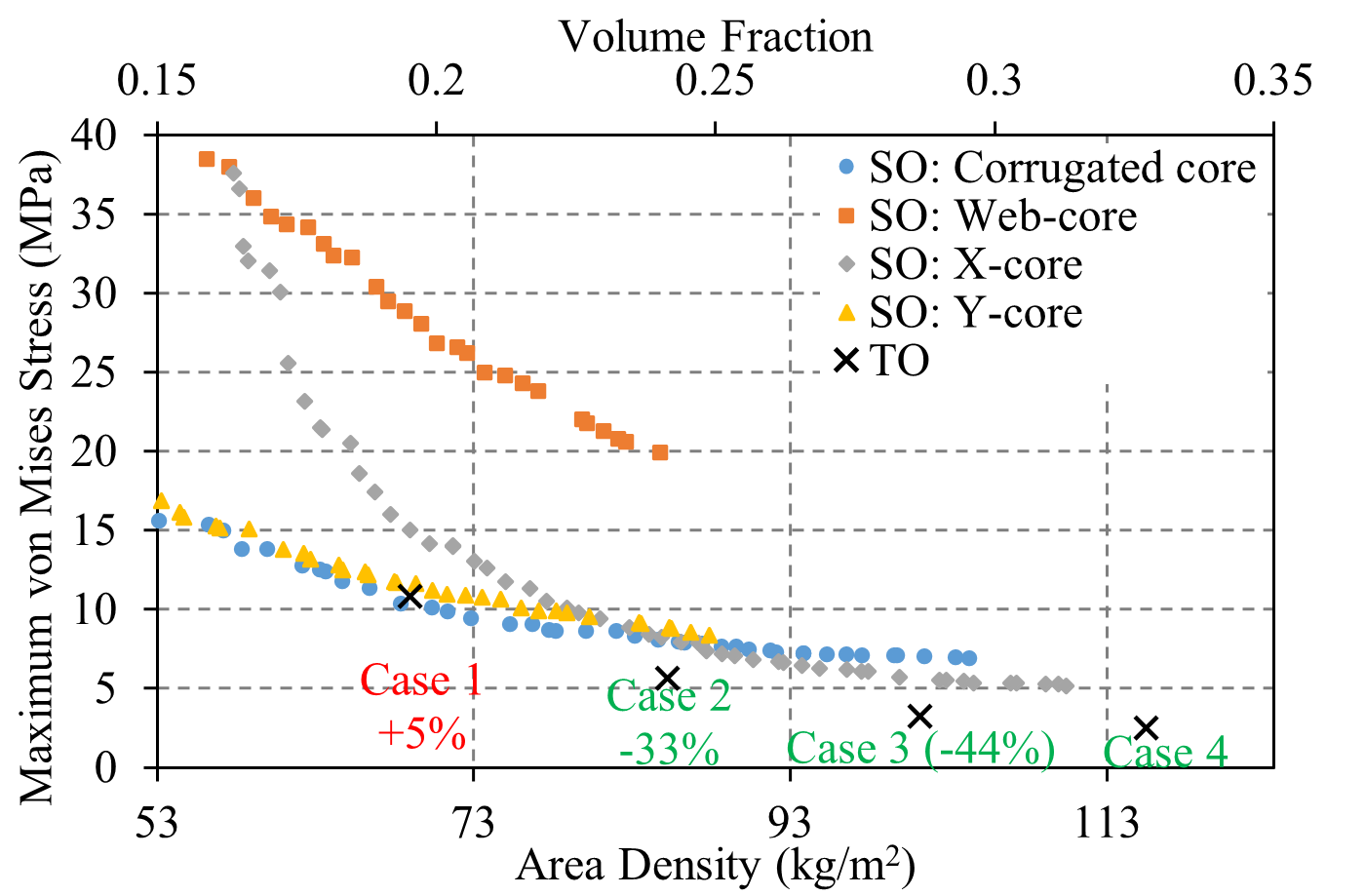}
    \caption{Subset of stress-based optimization results, highlighting TO solutions.}
    \label{noJstresszoomin}
\end{figure}

\begin{figure}[h]
    \centering
    \includegraphics[width=0.49\textwidth]{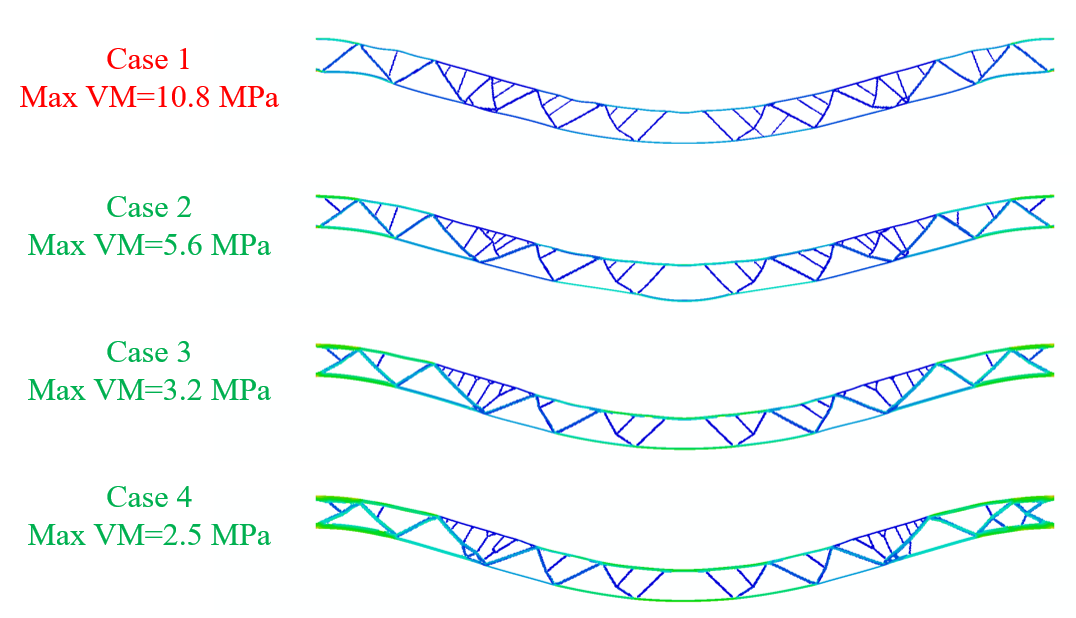}
    \caption{Topologies of optimized beams at different VF levels, shown in deformed state.}
    \label{noJabaqus}
\end{figure}

Figure \ref{noJstresszoomin} shows part of the results from Figure \ref{noJstress}, illustrating more closely the performance of TO solutions in comparison to the traditional sandwich beams. As can be seen, TO generates outstanding designs at higher area densities. When the volume fraction (VF) is below 0.2, TO is unable to improve the performance. When VF equals 0.2, TO exhibits 5\% higher stress than the best traditional sandwich beam with the same mass, which is in this case corrugated-core beam. Interestingly, increasing the VF to 0.24, the TO outperforms the best traditional-core beam by 33\%, which is quite significant. Increasing the VF further to 0.29, the improvement is 44\%. At VF of 0.33, the TO result is also outstanding, although beyond the area density range of the traditional beams. Figure \ref{noJabaqus} shows the topologies of these four cases. We can observe that they bear resemblance to corrugated-core beams, with some refinements in the form of additional trusses, primarily close to the upper surface where the load acts. With the increase in VF, thickness of the face plates at the supports and in the middle of the span increases in a gradual way, which are the places with large bending moments. At VF=0.2, thickness of the face plates is 5 mm at the edges and 4 mm in the middle of the span. At VF=0.33, face plates at edges are 13 mm thick and 6 mm in the middle. Plates in the core in all cases are at 2.5 mm and above. Furthermore, it could also be said that some TO results partially resemble X-core beams, most notably close to the midspan in case 1, and close to the right-hand edge in case 4. Thus, in a way, TO confirms that corrugated-core and X-core beams are optimal for stress minimization when VF ranges between 0.2 and 0.3. TO offered better performance of the beams, but with somewhat increased structural complexity than present in traditional prismatic sandwich beams. It is also interesting to observe in Figure \ref{noJabaqus} that although the loading and boundary conditions are fully symmetric from left to right, the topology has some minor asymmetries, most notable close to the supports in cases 1, 2 and 4. 

\begin{figure}[h]
	\centering
	\begin{subfigure}{0.99\linewidth}
		\includegraphics[width=\linewidth]{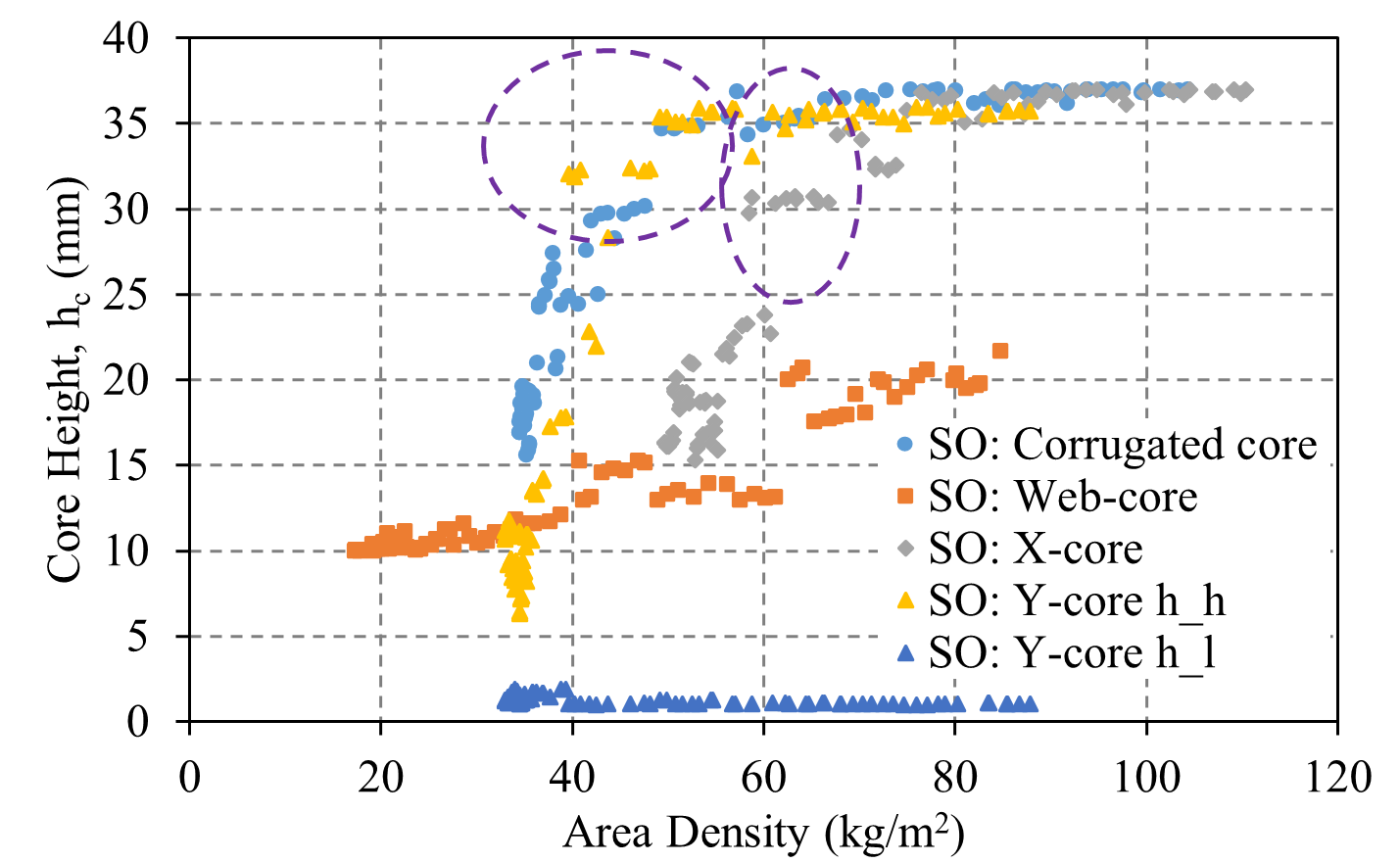}
		\caption{Core height.}
		\label{NoJStressHc}
	\end{subfigure}
 \vfill
	\begin{subfigure}{0.99\linewidth}
		\includegraphics[width=\linewidth]{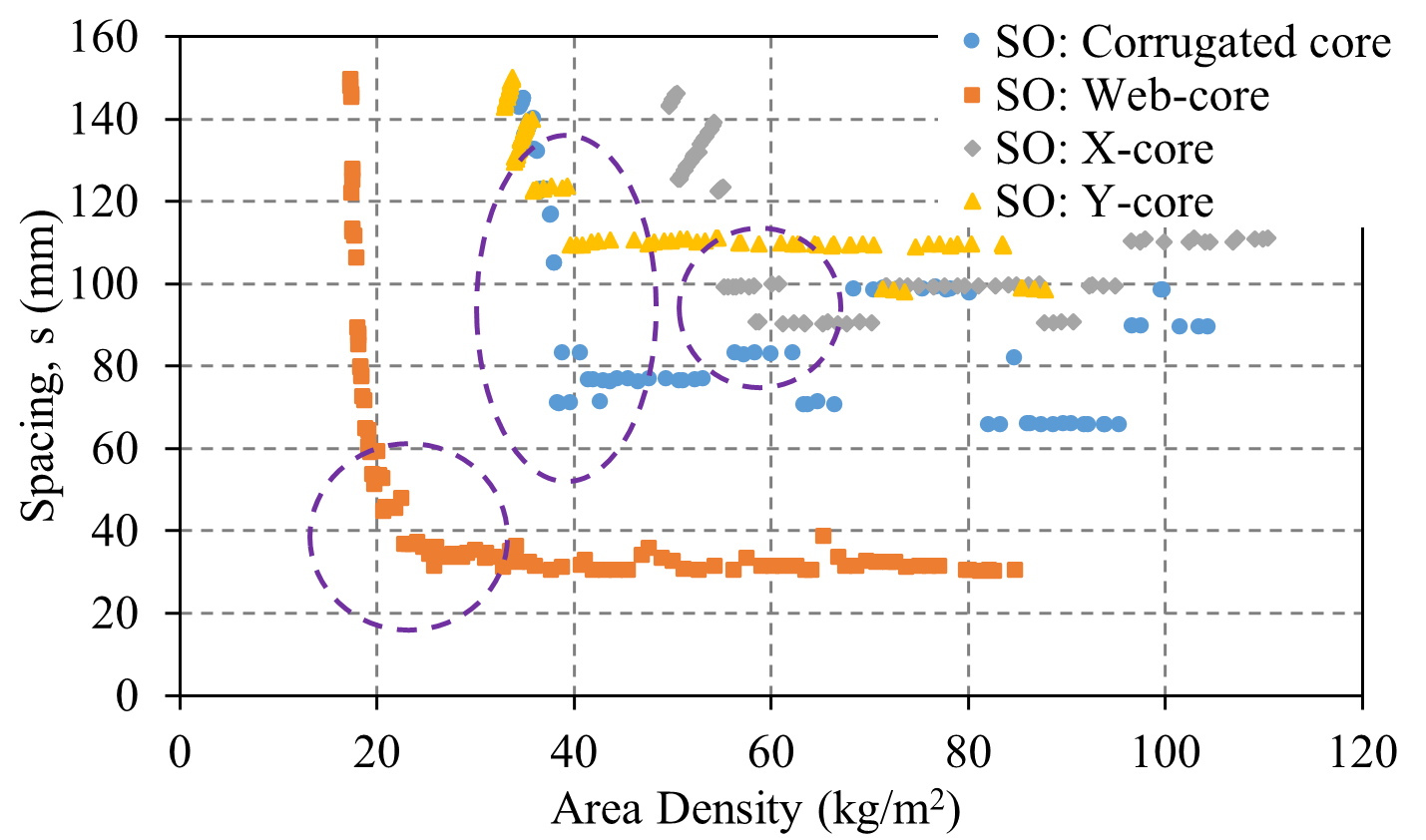}
		\caption{Spacing of profiles in the core.}
		\label{NoJStressS}
	\end{subfigure}
	\caption{Values of key design variables for traditional beams which exhibit knee points in stress-based optimization without joints.}
	\label{fig:NoJStressVariable}
\end{figure}

Arrows in Figure \ref{noJstress} highlight the ‘knee points’ in the Pareto fronts of the beams with traditional cores, with the exception of the web-core beam which has very uniform slope of the front. The knee points are caused by the variable bounds. Thickness of all the plates is at minimum for area densities below the knee points and gradually increases for heavier beams. On the other hand, core height is at maximum and spacing is at minimum for area densities above knee points, as seen in Figure \ref{fig:NoJStressVariable}, which results in heavier beams with lower maximum stress. As mentioned above, hl for the Y-core beam is at minimum for the entire Pareto front, forming effectively a corrugation, thus Y-core is not good for stress minimization. It is worth mentioning, however, that the benefits of Y-core sandwich structures are recognized elsewhere, featuring for example superior energy absorption capability in case of an impact when compared to corrugated-core sandwich panels \citep{st2015low} and X-core sandwich panels \citep{korgesaar2014steel}.

TO succeeds to generate better designs at higher area density levels. On the other hand, it struggles to discover beams lighter than 60 \(kg/m^3\) or those with a volume fraction (VF) of 0.2. With a VF of 0.23, TO exhibits a 5\% higher stress than the optimal nearby SO design point, but as the VF escalates to 0.28, the difference drastically plummets to -32\%. As can be observed in Figure \ref{noJabaqus}, the final topology significantly varies as VF ascends from 0.23 to 0.28, and thereafter evolves gradually with further increases in VF. This observation suggests that TO possesses an inherent threshold VF value, surpassing which significantly enhances its performance.

\subsubsection{Frequency optimization}
\begin{figure}[h]
    \centering
    \includegraphics[width=0.45\textwidth]{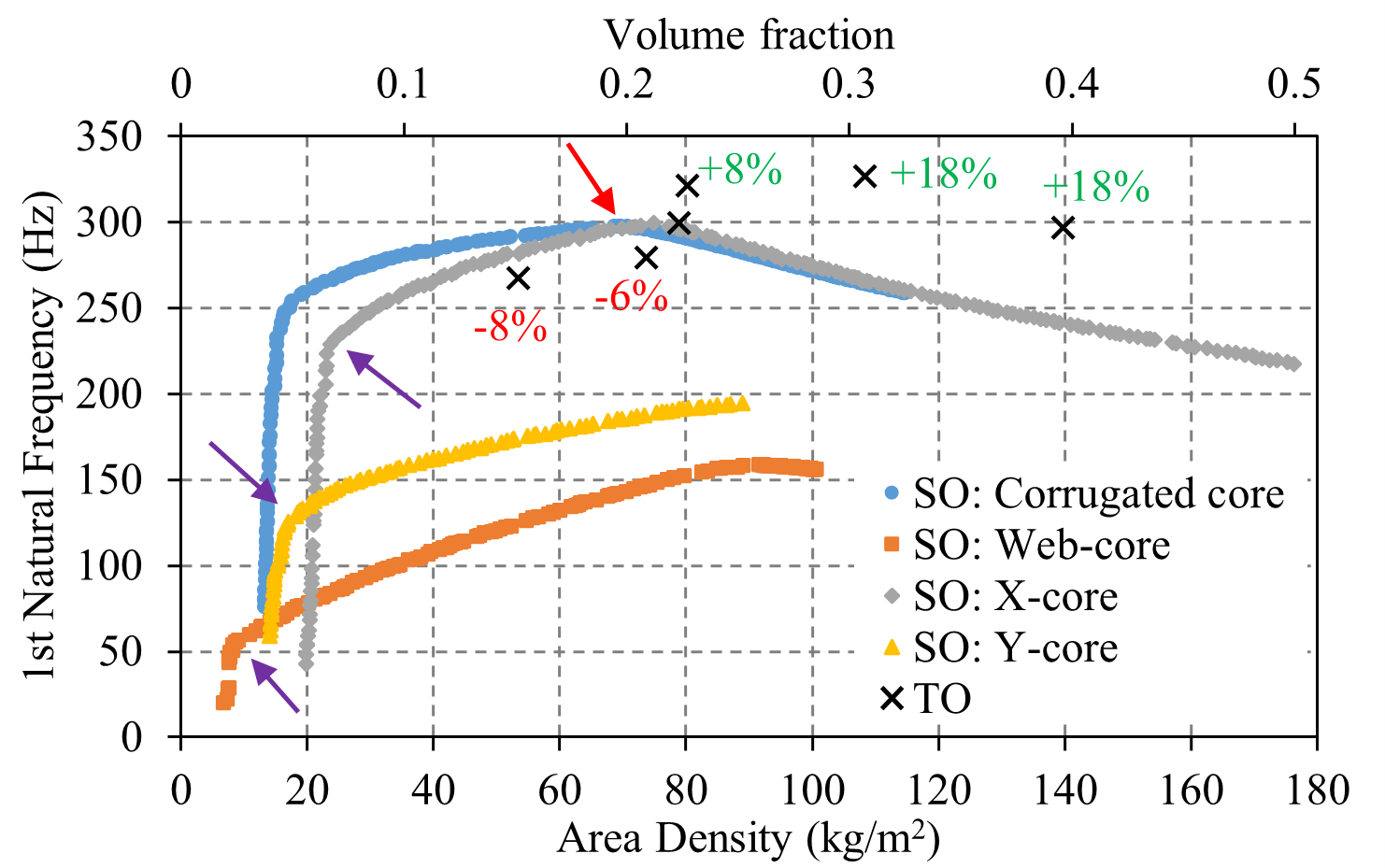}
    \caption{Optimal designs of frequency-based optimization without joints.}
    \label{freq1}
\end{figure}

To prevent resonance and reduce possibility of structural failures, the lowest natural frequency should be above excitation frequency. Figure \ref{freq1} shows the results of the first natural frequency maximization for traditional core sandwich beams and beams optimized via TO. Corrugated-core and X-core beams demonstrate similar performance, superior to the other two traditional core types. This superiority is attributed to the triangular shape of their cores, which ensures high transverse shear stiffness. While the upper part of the Y-core beams features triangular shape, the lower part is a vertical plate, which leads to lower transverse shear stiffness than in X-core and corrugated-core beams. Y-core beams in the Pareto front in Figure \ref{freq1} possess minimal height of the vertical plate, as in stress minimization case, suggesting this shape is not desirable. The web-core beam exhibits the worst performance in terms of frequency maximization, due to its lowest transverse shear stiffness among the four core types. Transverse shear stiffness was shown to strongly affect the response and strength of prismatic sandwich panels, and the web-core sandwich structures are particularly sensitive, see e.g., \cite{jelovica2012influence}
and \cite{jelovica2016eigenfrequency}.

\begin{figure*}[h]
  \includegraphics[width=\textwidth]{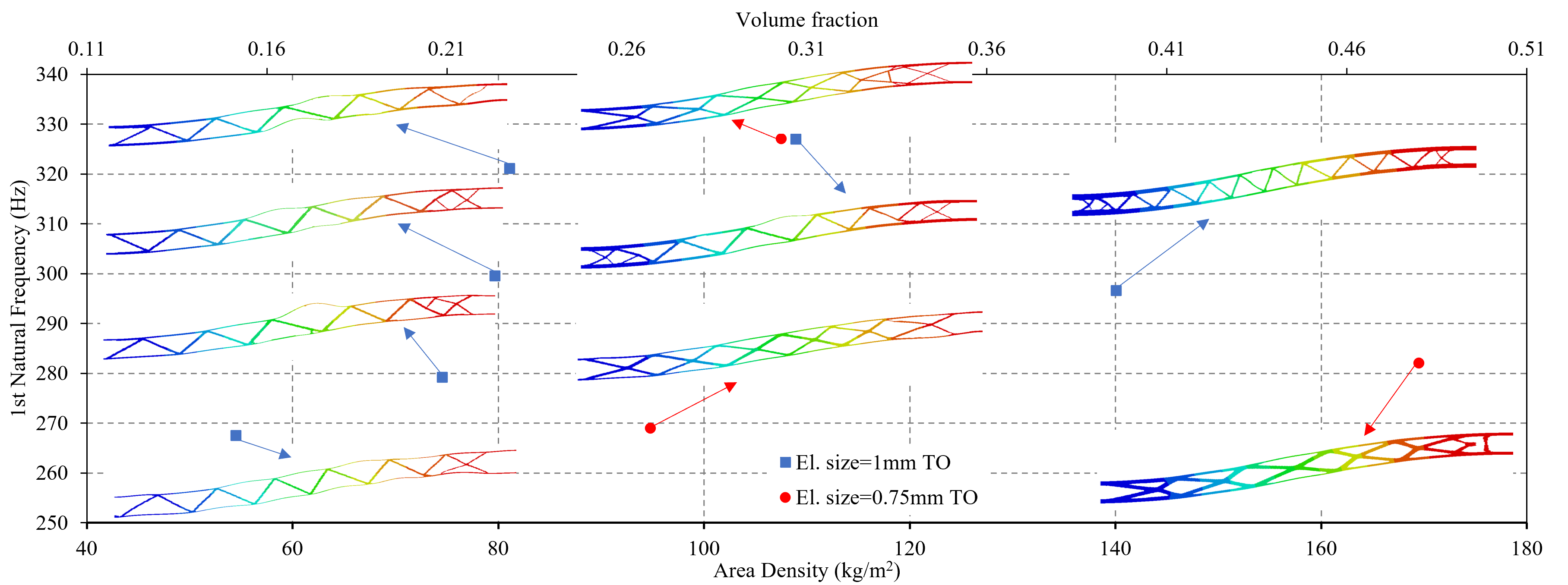}
  \caption{Frequency-based topology optimization of beams (color showing displacements in half of the structure for the first eigenmode).}
  \label{freq3}
\end{figure*}

\begin{figure}[h]
	\centering
	\begin{subfigure}{0.99\linewidth}
		\includegraphics[width=\linewidth]{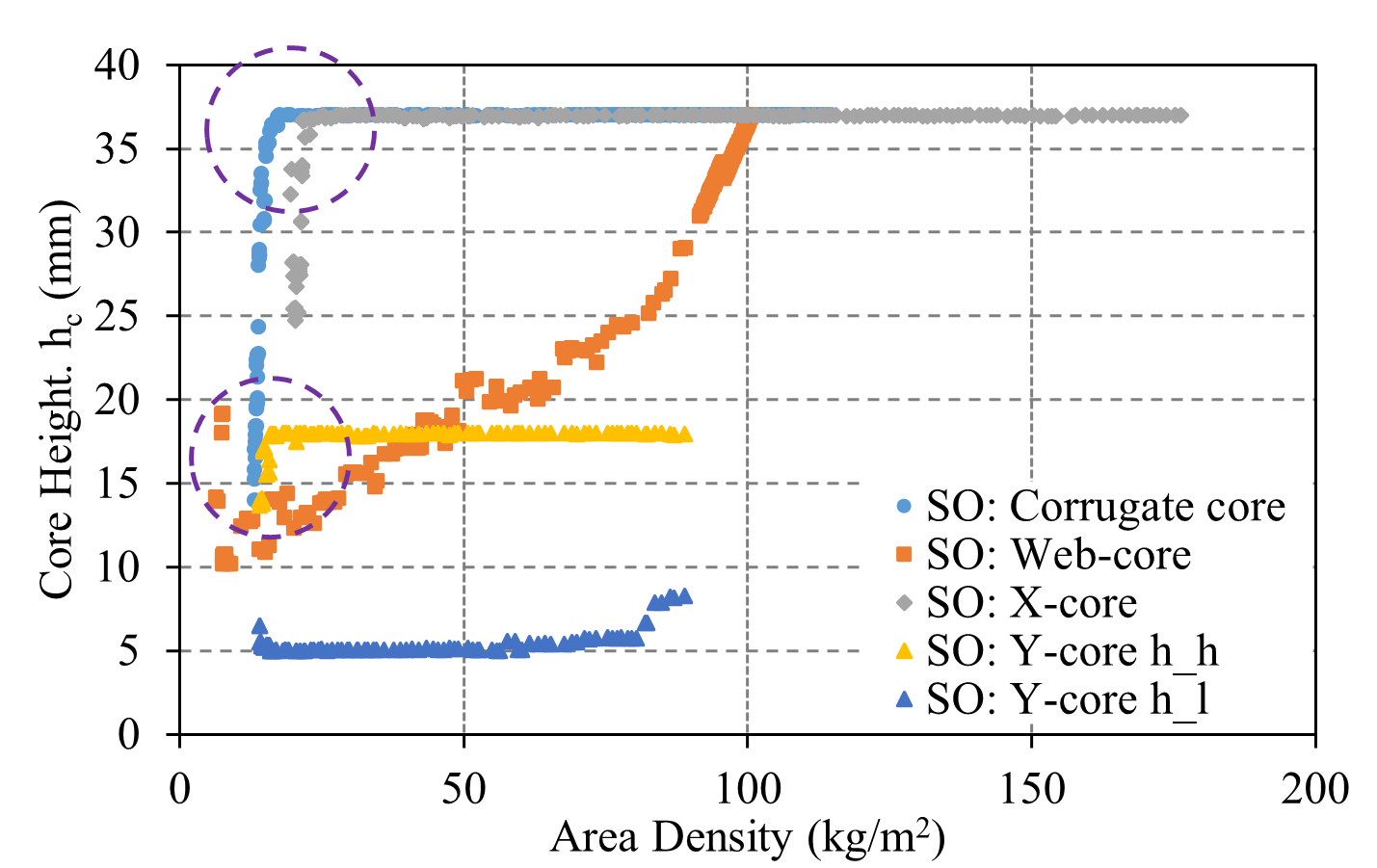}
		\caption{Core height.}
		\label{NoJStressHc}
	\end{subfigure}
 \vfill
	\begin{subfigure}{0.99\linewidth}
		\includegraphics[width=\linewidth]{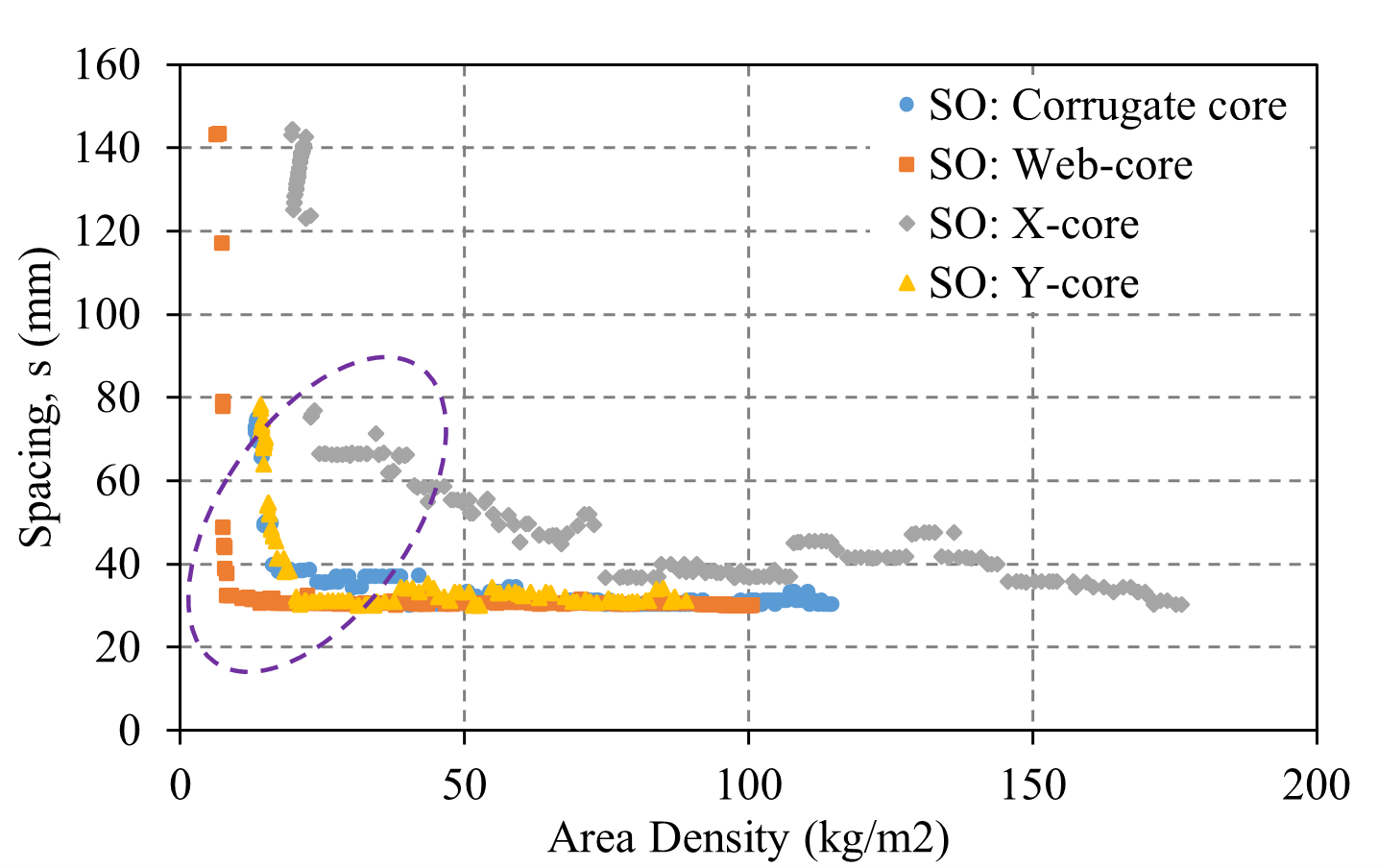}
		\caption{Spacing of profiles in the core.}
		\label{NoJStressS}
	\end{subfigure}
	\caption{Values of key design variables for traditional beams which exhibit knee points in frequency-optimized beams without joints.}
	\label{fig:FreqVariable}
\end{figure}

Figure \ref{freq1} showcases also the results of topology optimization in contrast to the traditional beams. For area densities exceeding 80 kg/m2, TO successfully identifies superior topologies, outperforming the traditional beams by 8-18\%. For area densities below 80 kg/m2 (VF less than 0.2), corrugated-core and X-core sandwich beams demonstrate superior performance. Notably, TO exhibits better performance when the volume fraction ranges between 0.23 and 0.4, similar to stress minimization in Section \ref{sec:stressnojoint}. Any deviations below or above this range tends to compromise the performance. Influence of element size on results is elaborated in Section \ref{subsec1}, where it is concluded that 1 mm is a good compromise between accuracy and computational expense, and is used throughout the study. Thus, the TO results in Figure \ref{freq1} are based on 1 mm element size. Figure \ref{freq3} shows those topologies, with a few additional ones that use an element size of 0.75 mm, for comparison. The smaller size (0.75 mm) grants TO a higher resolution, which in this case results in an X-shape core in a few instances, while the larger size (1 mm) generates predominantly corrugated-core beams. This confirms the high efficiency of both X-core and corrugated-core beams for frequency maximization and corresponds with the results for traditional prismatic sandwich beams. Further refinement of the mesh, however, does not significantly change the first natural frequency. In comparison to topologies based on stress minimization in Figure \ref{noJabaqus}, topologies in Figure \ref{freq3} feature simpler shapes without miniscule struts, which are unnecessary in this case due to the lack of localized loading. For the 1 mm element size, the lowest face plate thickness in TO is 1 mm (as set through active elements), while the plates in the core are on average 2.5 mm thick, although there are localized parts where that thickness drops to 1.5 mm, which is lower than allowed for traditional beams. However, as it occurs only locally in few cases, the thickness bounds of the design space are preserved overall. 
Similar to Section \ref{sec:stressnojoint}, ‘knee points’ are observed in Pareto fronts of the traditional beams in Figure \ref{freq1}, indicated with arrows; however, in case of frequency optimization there are now two sets of knee points per front. As the area density increases, the first knee point is caused by the fact that the stiffness of the beams cannot increase anymore at the same rate, since the core height has reached its maximum and spacing minimum, as seen in Figure \ref{fig:FreqVariable}. Increase in stiffness from that point on is due to the increase in thickness of all the plates. Finally, going beyond the second knee point, the increase in mass is more significant than increase in stiffness, which ultimately leads to the decrease of natural frequency. It is interesting to note in Figure \ref{fig:FreqVariable} that the core height of the web-core beam continues to gradually increase through the entire area density range, but the web-plate spacing drops to the minimum at the beginning, causing pronounced first knee point. 

Frequency optimization of the traditional beams is carried out using shell elements (S8R). To validate the results, selected designs are analyzed using solid elements (C3D8R), and the result is shown in Figure 15, demonstrating an excellent match. Another validation is performed by comparing the FEA solver in MATLAB (TO) and Abaqus (SO). The results are given in Table \ref{tab2}. The difference in frequency is within 2\%, therefore the accuracy of the MATLAB solver is sufficient for predicting the natural frequency.

\begin{figure}[h]
    \centering
    \includegraphics[width=0.45\textwidth]{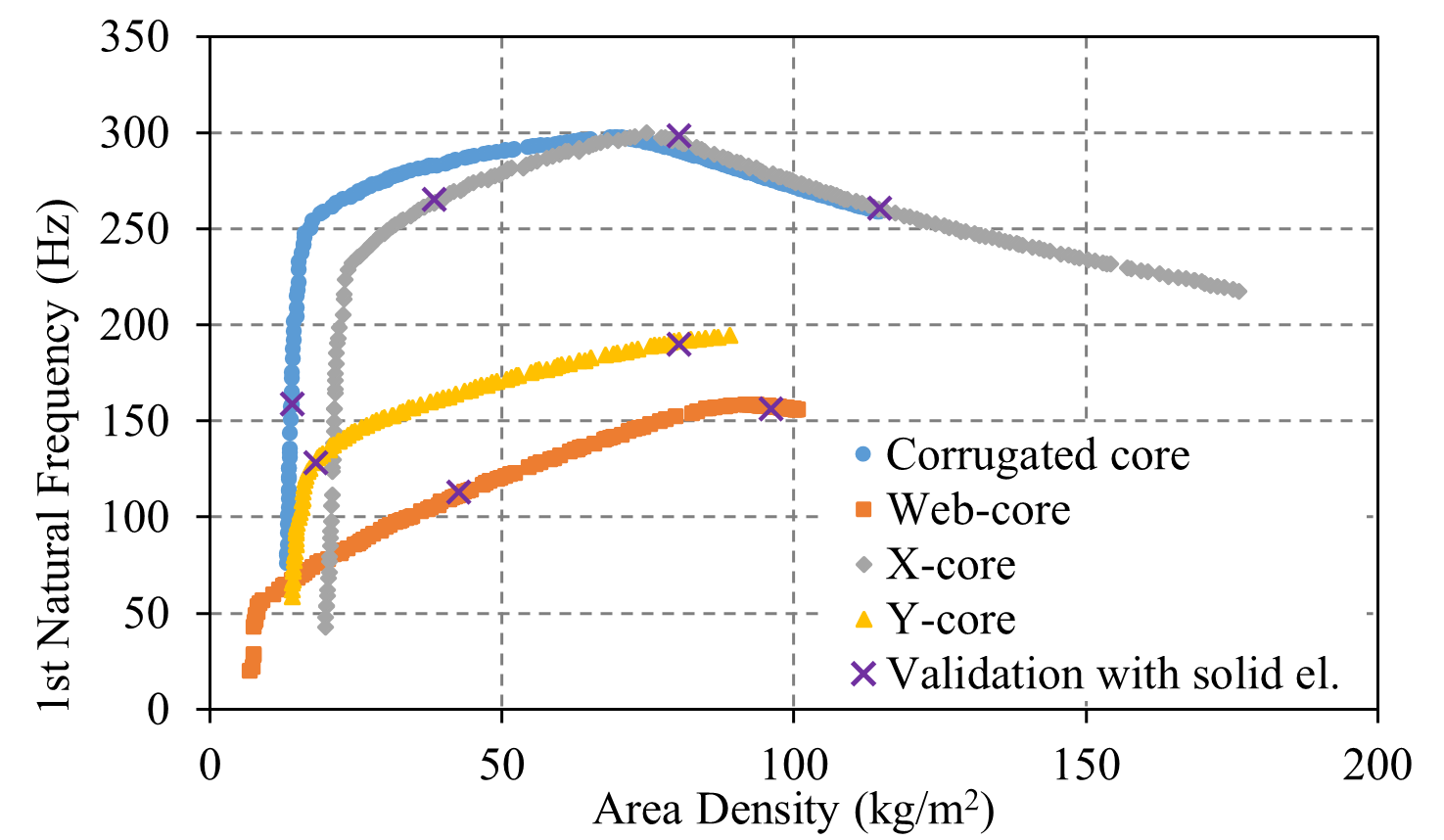}
    \caption{Validation of SO results for frequency-optimized beams without joints.}
    \label{Freqvali}
\end{figure}

\begin{table}[h]
\begin{center}
\begin{minipage}{325pt}
\caption{Comparison of FEA solvers in TO and SO.}\label{tab2}
\begin{tabular}{@{}llll@{}}
\toprule
El. size [mm] & f(Abaqus)[Hz] & f(MATLAB)[Hz] & Diff. \\
\midrule
1.5    & 296  & 296  & 0.18\% \\
1.0    & 326   & 323  & 0.92\% \\
0.75    & 326   & 321 & 1.60\% \\
0.6    & 326   & 320  & 1.84\% \\
0.5    & 327   & 321  & 1.83\%  \\
\bottomrule
\end{tabular}
\end{minipage}
\end{center}
\end{table}

\subsection{Stress optimization on sandwich beams with joints}\label{sec:joint}
\begin{figure}[h]
    \centering
    \includegraphics[width=0.45\textwidth]{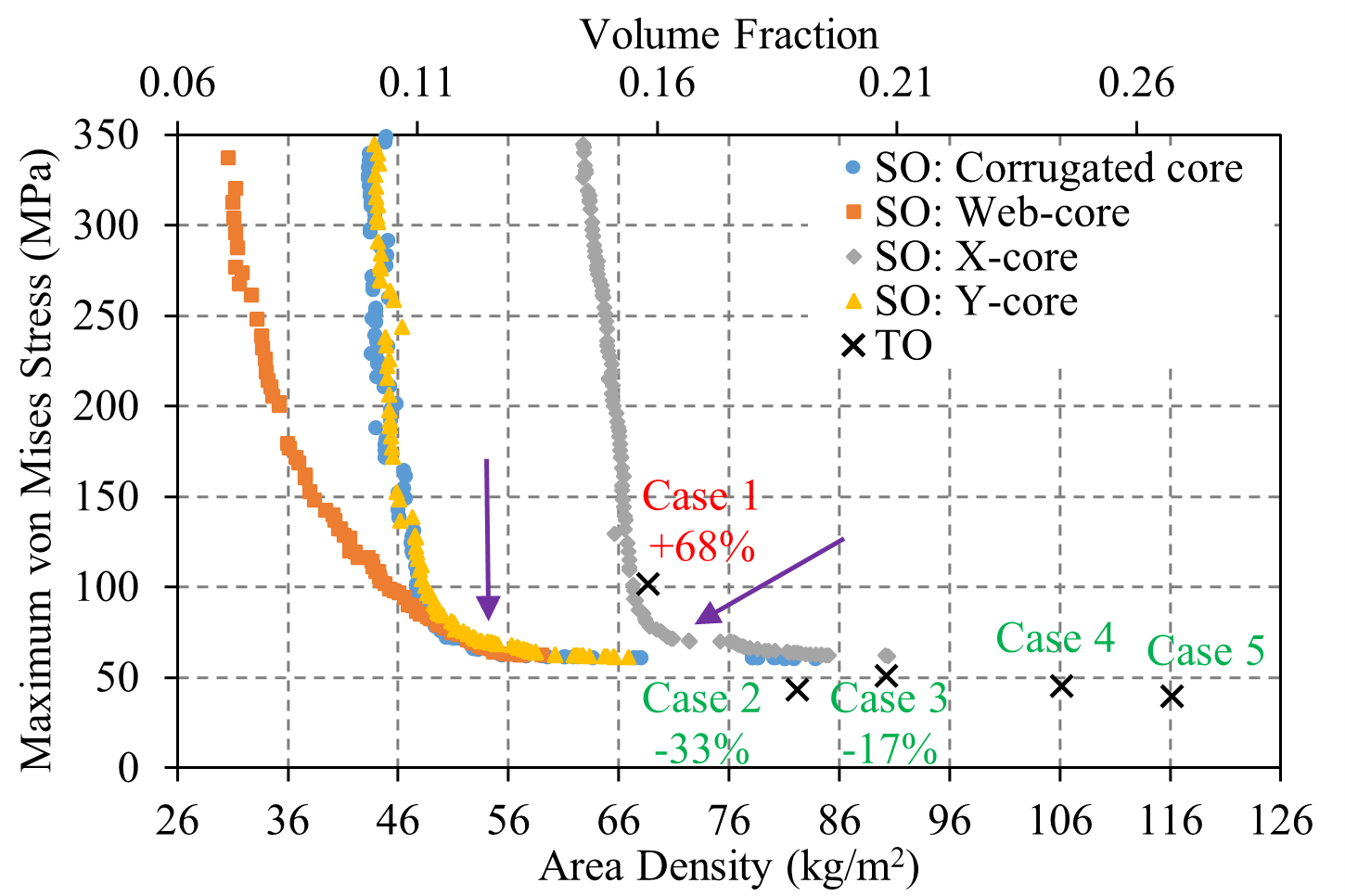}
    \caption{Stress optimization results for the sandwich beams with joints.}
    \label{bond1}
\end{figure}

\begin{figure}[h]
    \centering
    \includegraphics[width=0.49\textwidth]{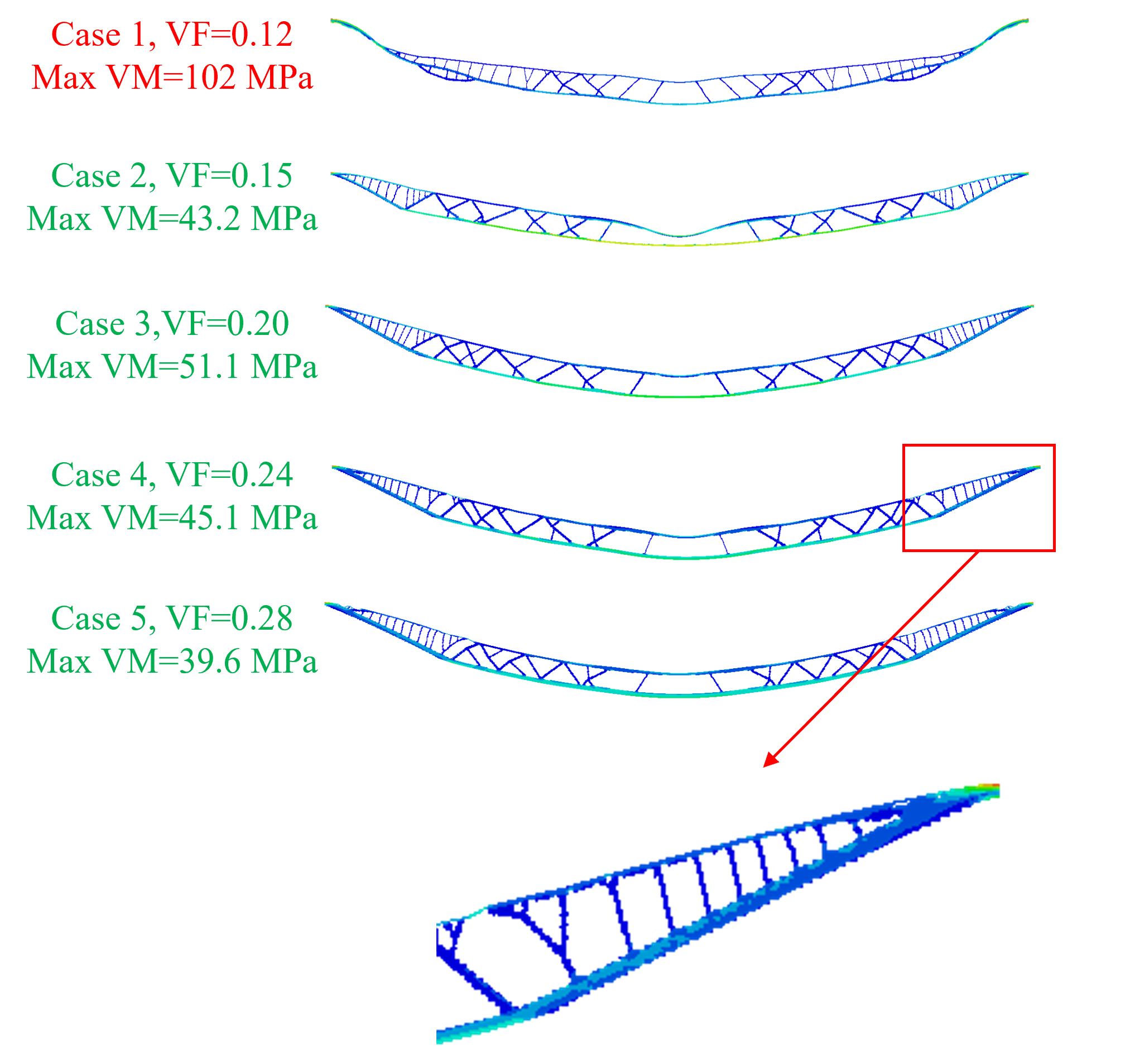}
    \caption{Topologically optimized beam structures with joints at different VF levels in deformed state.}
    \label{Jabaqus}
\end{figure}

Figure \ref{bond1} shows the results of the stress optimization of beams with joints. A five-fold increase in stress can be observed when compared to beams without joints (Figure \ref{noJstress}). This is due to the high localized stresses in the joints, in both the conventional joint used for traditional core beams and innovative joints created through topology optimization. The Pareto fronts for the Y-core and corrugated-core beams show high degree of overlap, for the same reason explained in Section \ref{sec:stressnojoint}. At higher mass levels, TO outperforms optimized traditional beams by 33\%. At the lower end of area densities (cases 1), TO experiences 68\% higher stress than the best traditional sandwich beam with the same mass, which is the corrugated core beam in this case. Thus, the TO does not necessarily converge to the global optimum, which is not surprising given that local optimizer (MMA) is used for computational efficiency, and the targeted VF is 0.1 in Case 1, which leaves a very small amount of material for TO to distribute. As can be observed in Figure \ref{Jabaqus}, the overall layout of the final topology significantly alters to another shape with triangular-shaped joint sections at the ends from Case 2 as VF increases. Case 1 is considered an outlier in this case. Global optimizers could potentially find even better designs, but they are often computationally prohibitive within the TO framework for a design domain with such a large aspect ratio.

Topologies of beams optimized through TO are different when joints are considered, as can be seen when comparing Figure \ref{Jabaqus} with Figure \ref{noJabaqus}. Without the joints the core of the TO beams resembles predominantly corrugated-core and only to small extent X-core beams, while with joints the situation is reversed. In addition, with joints present, parts of the beams feature vertical struts, resembling web-core beams. It is interesting to note that the generated topology of the joints is very simple for the low VFs, and as mentioned above, is not globally optimal. At higher VFs vertical struts are added to support the top plate of the joint, which is effective in decreasing stress, but results in a joint that is more difficult to manufacture than the joint previously proposed in the literature. The ‘traditional’ joint has span (in x direction) of 120 mm; the TO joint spans from 130 mm (Case 2) to 206 mm as the volume fraction increases. The thickness of the plates in the traditional joint is up to 3 mm; TO joint is up to 6 mm (Case 4 and 5). In that case, adjacent to the edge the top face is 2 mm thick while vertical struts are between 2 mm and 3 mm thick.
 
One interesting finding is that TO was unable to generate optimal solutions at all area density levels. Sub-optimal solutions have been generated for low-to-intermediate area density levels, and in some cases, they even need manual modification. This issue is considered to originate from the relatively large aspect ratios of the design domains (22.2 in Section \ref{sec:nojoint} and 27.6 in Section \ref{sec:joint}). The aspect ratio used for TO in the literature is usually within one and three \citep{liu2014efficient},\citep{deng2021efficient}. We found that certain optimizers are sensitive to the aspect ratio, e.g., OC was unable to converge to a solution for cases in Section \ref{sec:joint}. The aspect ratio influences the value distribution of the sensitivity information which is crucial for gradient-based optimizers like OC and MMA. Although this TO framework did not always converge to global optima, it offers substantial improvements while confirming that the corrugated core and X-core beams are very good core types.

\begin{figure}[h]
	\centering
	\begin{subfigure}{0.99\linewidth}
		\includegraphics[width=\linewidth]{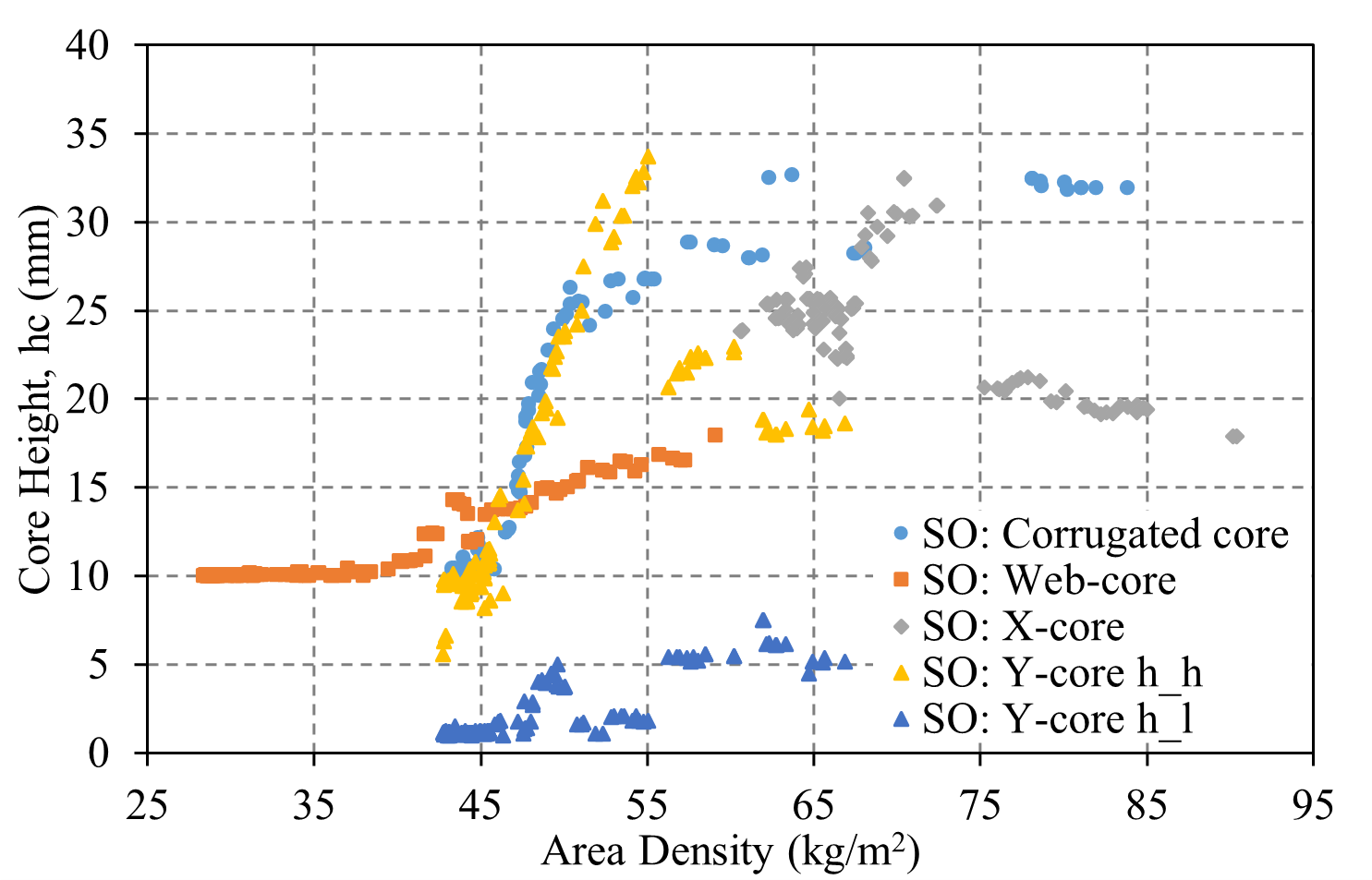}
		\caption{Core height.}
		\label{JStressHc}
	\end{subfigure}
 \vfill
	\begin{subfigure}{0.99\linewidth}
		\includegraphics[width=\linewidth]{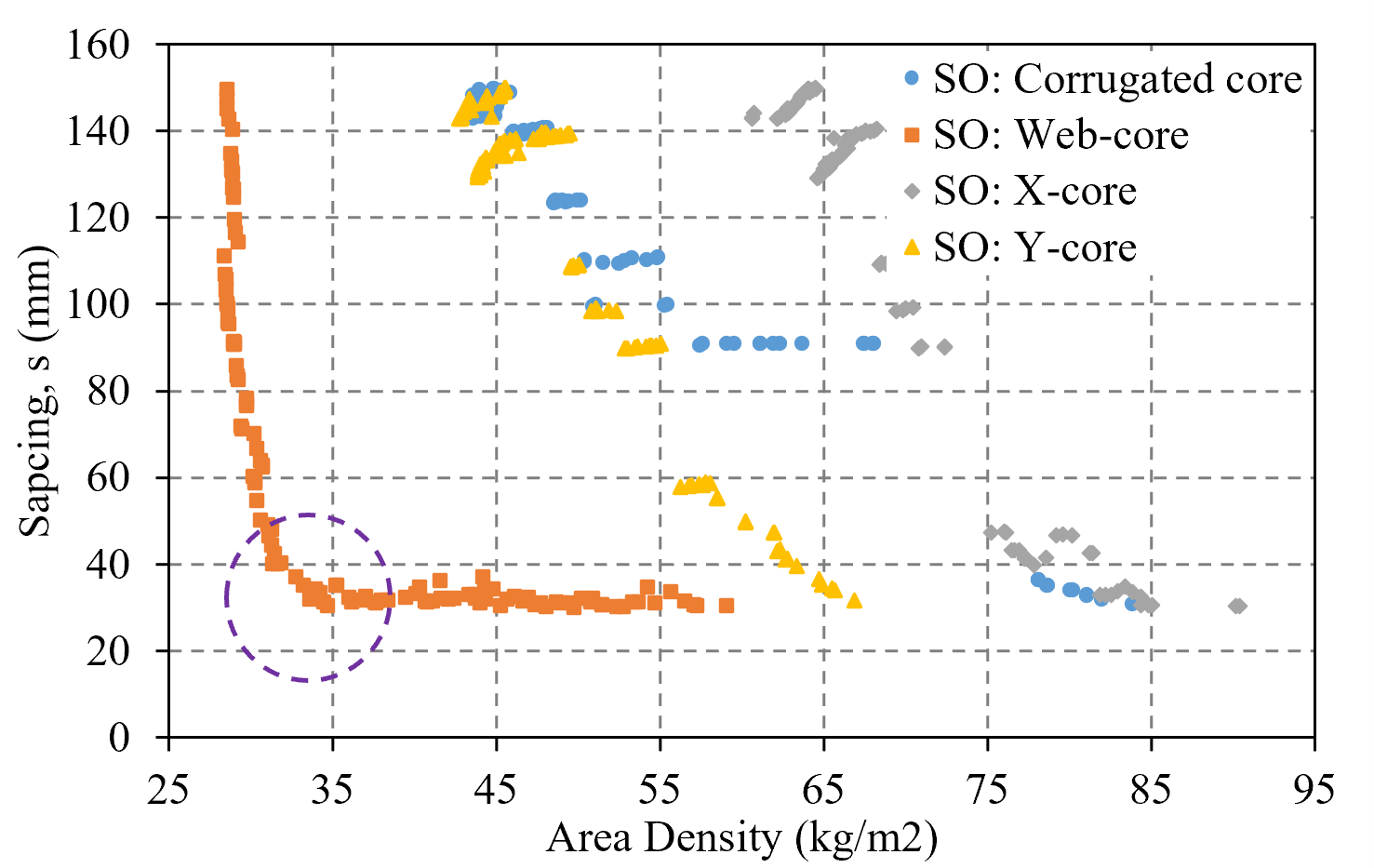}
		\caption{Spacing of profiles in the core.}
		\label{JStressS}
	\end{subfigure}
	\caption{Values of key design variables for traditional beams which exhibit knee points in stress-optimized beams with joints.}
	\label{fig:JStressVariable}
\end{figure}

Figure \ref{fig:JStressVariable} shows core height and core spacing of the traditional beams with joints. We can see that in this case the change in variables is somewhat different than the case without joints (Figure \ref{fig:NoJStressVariable}). With joints, the solutions in the Pareto fronts for corrugated-, X- and Y-core beams feature lower core height at higher area densities than without joints. Similarly, with joints, these beam types have lower spacing between profiles in the core, in comparison to the case without joints. One possible explanation is that lower core height effectively reduces stresses at the joint. Knee points are still present at the Pareto fronts (Figure \ref{bond1}), which are characterized by a peak or extreme in variable values (Figure \ref{fig:JStressVariable}). As above, web-core beam shows a gradual change in variable values. 

\subsection{Influence of mesh density}\label{subsec1}
Attaining mesh independence in topology optimization is theoretically achievable through the incorporation of an effective filter, such as the density filter employed in this research \citep{svanberg2013density}. However, in this study, perfect mesh independence has not been realized in the context of frequency optimization and stress optimization, with the topology showing qualitative differences between different element sizes. This phenomenon could be attributed to the fixed setting of the filter radius within the density filter, which in itself can act as a control parameter for the final topology \citep{yan2022topology}. The filter radius needs to be set to a fixed percentage of the characteristic length of the design domain to guarantee mesh independence. As per \cite{sigmund200199}, TO ideally does not show qualitative differences in topologies between different element sizes. However, there could be benefits even if the solution is not entirely mesh-independent, since in such case better novel solutions could be obtained. Consequently, the impact of mesh density on the results warrants further explanation.

\begin{figure}[H]
    \centering
    \includegraphics[width=.45\textwidth]{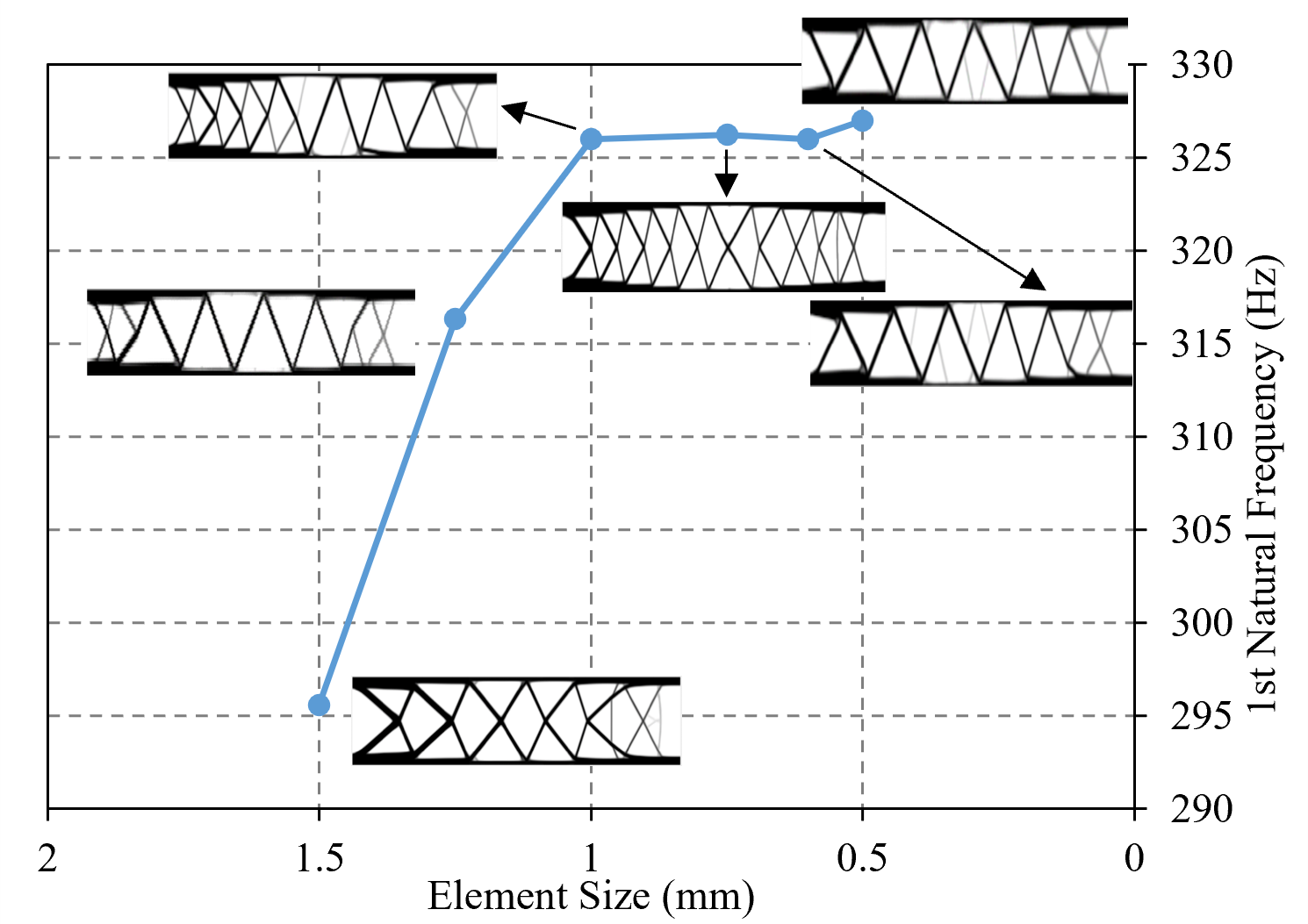}
    \caption{Mesh sensitivity analysis in frequency-based TO with VF=0.3.}
    \label{sens}
\end{figure}

In order to achieve a balance between computational expense and accuracy, TO has been executed with varying mesh sizes, and the resultant optimal topologies and natural frequencies are illustrated in Figure \ref{sens}. The filter radius is set to a relatively low value (r = 2) to facilitate generation of more intricate solutions as the element size diminishes. As the figure shows, element size between 1 mm and 0.5 mm results in basically constant first natural frequency. Despite the expectation that the smallest element size would allow the most complex topology, this is actually observed for element length of 0.75 mm, and not in the case of smaller elements. However, as the element size decreases, dimensions of the matrices escalate dramatically, thereby incurring a significant computational cost. In an attempt to achieve the balance between accuracy and computational cost, element size of 1 mm is used in the rest of the study, except for the frequency optimization where 0.75 mm element is used for comparison, as shown in Figure \ref{freq3}. The filter radius is kept at r = 3 throughout the study.

\subsection{Computational cost} \label{comp}
The computational expense of both SO and TO for frequency and stress optimization (without joints) has been summarized in Table \ref{computation}. The SO utilizes Abaqus R2018 as a solver for estimating the values of stress and frequency. This process necessitates an average of 7 seconds, which comprises of 4 seconds dedicated to the activation of the solver, 2 seconds devoted to the execution of the FEA, and the final second for post-processing. The frequency-based TO features a full matrix implementation and the implicitly restarted Lanczos method as the solver. Notably, due to the significant size of the characteristic matrix, it requires an average of 2000 seconds for each FEA run. On a positive note, the objective function converges more readily, necessitating merely 300 iterations. On the contrary, stress-based TO utilizes a sparse matrix implementation and a similar high-performance sparse solver to expedite the computation \citep{stewart2002krylov}. This results in a reduction of FEA time to only 2.3 seconds in TO, which is approximately a thousand times faster than in the case of the dense matrix implementation in frequency-based TO. Nevertheless, the sensitivity information obtained from stress-based TO is less stable compared to the frequency-based TO, hence necessitating a larger number of iterations for convergence (more than 2000). Despite this, the computation still benefits from the advancements in sparse matrix computation, and TO utilizes half of the time compared to SO in this particular scenario because SO needs at least 12,000 designs per algorithm run, for each traditional beam type.

\begin{table}[h]
\begin{center}
\begin{minipage}{170pt}
\caption{Computational cost of the TO and SO.}\label{computation}
\begin{tabular}{@{}lllll@{}}
\toprule
& \multicolumn{2}{@{}c@{}}{Frequency-based} & \multicolumn{2}{@{}c@{}}{Stress-based} \\\cmidrule(lr){2-3}\cmidrule(lr){4-5}%
~ & SO & TO & SO & TO \\
\midrule
No. of FEA    & 12,000   & 300  & 12,000 & 30,000  \\
Times per FEA [s]    & 7   & 2,000  & 7 & 2.3  \\
Total time [min]    & 1,400   & 10,000  & 1,400 & 1,150  \\
Difference    & ~   & +614\%  & ~ & -18\%  \\
\bottomrule
\end{tabular}
\end{minipage}
\end{center}
\end{table}

These findings suggest that TO can potentially be faster than SO, even though it is conventionally deemed to be more computationally demanding than sizing optimization \citep{holmberg2013stress}. The acceleration of TO has been accomplished through the introduction of advanced FEA solvers and sparse matrix implementation in the stress-based TO in this study. Further acceleration of TO can be achieved by parallel computing \citep{aage2015topology}, neural networks (NN) \citep{chandrasekhar2021tounn} and non-gradient-based optimizers \citep{biyikli2015proportional}.

\section{Conclusions}\label{sec5}   
A few prominent types of traditional prismatic sandwich beams are optimized in this study, namely X-, Y-, corrugated- and web-core beams. Von Mises stress in beams under pressure loading is minimized and natural frequency is maximized at different mass levels. The performance and characteristics of the optimal beams are compared for the first time with the beams created through topology optimization (TO). In addition, stress is minimized in beams with joints (or connections) to the side supports. The considered traditional beams feature a joint proposed in literature, whose thickness is optimized but topology remains unchanged. This is compared to the case where topology optimization identifies the joint topology together with that of the beam. This is the first study that considers the response of metallic sandwich structures with joints, and seeks better performance through topology optimization framework. The literature is extremely scarce when it comes to design of joints for metallic sandwich structures, which is a big gap that needs to be addressed to help in allowing the use of sandwich panels in real-life structures. Efficiency and computational cost of the two optimization approaches are discussed. Sizing and topology optimization are performed using NSGA-II and MMA algorithms, respectively. Through the use of active elements and filter radius, the minimum thickness of faces and struts in TO is principally equivalent to sizing optimization, although somewhat more complex topologies of the core are produced, which needs to be considered in conjunction with manufacturing capabilities in future work.  

This comprehensive analysis contributes the following original insights:

\begin{enumerate}
\item Structural performance of the beams is improved by up to 44\% through topology optimization but at the cost of increased topological complexity which would require advanced manufacturing
\item Beams created through topology optimization outperform the traditional beams for volume fractions between 0.2 and 0.4; outside that range the traditional beams are better
\item Topologies of the newly created beams confirm that the traditional corrugated-core and X-core sandwich structures are in fact outstanding; they bear many resemblances to the topologies of the new cores
\item Joints of the sandwich beams are critical for performance improvements and need to be considered in optimization; TO offered improved joints with relatively simple topologies
\item Optimized traditional beams generally demonstrate excellent performance considering they can be manufactured in relatively simple way; corrugated-core and X-core beams are the best for volume fractions between 0.15 and 0.4
\item Control parameters, such as element size and filter radius, significantly impact the performance of TO
\end{enumerate}

\section{CRediT authorship contribution statement}
\textbf{Shengyu Yan}: Conceptualization, Methodology, Software, Model generation, Formal Analysis, Writing - original draft, Writing - review \& editing. \textbf{Jasmin Jelovica}: Supervision, Conceptualization, Methodology, Funding acquisition, Project administration, Validation, Writing - original draft, Writing - review \& editing.

\section{Acknowledgement}
The authors acknowledge the financial support by Natural Sciences and Engineering Research Council of Canada (NSERC) [grant numbers IRCPJ 550069-19 and RGPIN-2017-04509]. This research was supported in part by computational resources and services provided by Advanced Research Computing at The University of British Columbia.

\section{Declarations}
\textbf{Conflict of interest}: The authors declare that they have no conflict of interest.\\ \textbf{Replication of results}: The code and data are available from the authors upon reasonable request.








\bibliographystyle{unsrtnat}
\bibliography{snbibliography}

\end{document}